\documentclass{article}
\usepackage{latexsym}
\usepackage{amssymb}
\parskip=\medskipamount
\newtheorem{definition}{Definition}[section]
\newtheorem{lemma}[definition]{Lemma}
\newtheorem{theorem}[definition]{Theorem}
\newtheorem{proposition}[definition]{Proposition}
\newtheorem{corollary}[definition]{Corollary}
\newtheorem{remark}[definition]{Remark}

 1   
 1  

\def\QED{\hskip0.1em\hfill\null\ \null\nobreak\hfill
\kern3pt\lower1.8pt\vbox{\hrule\hbox
{\vrule\kern1pt\vbox{\kern1.7pt \hbox{$\scriptstyle
QED$}\kern0.2pt}\kern1pt\vrule}\hrule}}

\def\lcf{\lbrack\! \lbrack}
\def\rcf{\rbrack\! \rbrack}

\mathchardef\za="710B  
\mathchardef\zb="710C  
\mathchardef\zg="710D  
\mathchardef\zd="710E  
\mathchardef\zve="710F 
\mathchardef\zz="7110  
\mathchardef\zh="7111  
\mathchardef\zvy="7112 
\mathchardef\zi="7113  
\mathchardef\zk="7114  
\mathchardef\zl="7115  
\mathchardef\zm="7116  
\mathchardef\zn="7117  
\mathchardef\zx="7118  
\mathchardef\zp="7119  
\mathchardef\zr="711A  
\mathchardef\zs="711B  
\mathchardef\zt="711C  
\mathchardef\zu="711D  
\mathchardef\zvf="711E 
\mathchardef\zq="711F  
\mathchardef\zc="7120  
\mathchardef\zw="7121  
\mathchardef\ze="7122  
\mathchardef\zy="7123  
\mathchardef\zf="7124  
\mathchardef\zvr="7125 
\mathchardef\zvs="7126 
\mathchardef\zf="7127  
\mathchardef\zG="7000  
\mathchardef\zD="7001  
\mathchardef\zY="7002  
\mathchardef\zL="7003  
\mathchardef\zX="7004  
\mathchardef\zP="7005  
\mathchardef\zS="7006  
\mathchardef\zU="7007  
\mathchardef\zF="7008  
\mathchardef\zW="700A  

\newcommand{\be}{\begin{equation}}
\newcommand{\ee}{\end{equation}}

\newcommand{\bea}{\begin{eqnarray}}
\newcommand{\eea}{\end{eqnarray}}
\newcommand{\beas}{\begin{eqnarray*}}
\newcommand{\eeas}{\end{eqnarray*}}
\newcommand{\Z}{{\mathbb Z}}
\newcommand{\R}{{\mathbb R}}

\newcommand{\1}{{\mathbf 1}}

\newcommand{\D}{{d}}

\newcommand{\we}{\wedge}
\newcommand{\nn}{\nonumber}

\newcommand{\pa}{\partial}
\newcommand{\ti}{\times}

\newcommand{\X}{{\mathfrak X}}

\newcommand{\Ll}{{\pounds}}
\def\lna{\lbrack\! \lbrack}
\def\rna{\rbrack\! \rbrack}

\def\zT{{\cal T}}

\begin{document}

\title{Jacobi structures on affine bundles}
\author{J. Grabowski$^{1}$, D. Iglesias$^{2}$,  J.C.
Marrero$^{3},$ E. Padr\'on$^{3}, $ P. Urba\'nski$^{4}$
\\[10pt]
 {\small\it$^1$ Mathematical Institute, Polish Academy of
Sciences}\\[-6pt] {\small\it \'Sniadeckich 8, P.O.Box 21, 00-956
Warsaw, Poland }\\[-6pt]{\small\it E-mail:
jagrab@impan.gov.pl}\\[-4pt] {\small\it $^2$Instituto de
Matem\'aticas y F\'\i sica Fundamental} \\[-6pt] {\small\it Consejo
Superior de Investigaciones Cient\'\i ficas} \\[-6pt] {\small\it Serrano 123, 28006
Madrid, SPAIN}
\\[-6pt] {\small\it E-mail:
iglesias@math.psu.edu}\\[-4pt]{\small\it $^3$Departamento de
Matem\'atica Fundamental, Facultad de Matem\'aticas}\\[-6pt]
{\small\it Universidad de la Laguna, La Laguna} \\[-6pt]
{\small\it Tenerife, Canary Islands, SPAIN}\\[-6pt] {\small\it
E-mail: jcmarrer@ull.es, mepadron@ull.es }\\[-4pt]
{\small\it$^4$Division of Mathematical Methods in Physics,
University of  Warsaw}\\[-6pt] {\small\it Ho\.za 74, 00-682 Warsaw,
Poland }\\[-6pt] {\small\it E-mail: urbanski@fuw.edu.pl}}
\date{\empty}

\maketitle
\begin{abstract}
We study affine Jacobi structures (brackets) on an  affine bundle
$\pi:A\to M$, i.e. Jacobi brackets that close on affine functions.
We prove that if the rank of $A$ is non-zero, there is a
one-to-one correspondence between affine Jacobi structures on $A$
and Lie algebroid structures on the vector bundle
$A^+=\bigcup_{p\in M}Aff(A_p,\R)$ of affine functionals. In the
case $rank\; A=0$, it is shown that there is a one-to-one
correspondence between affine Jacobi structures on $A$ and local
Lie algebras on $A^+$. Some examples and applications, also for
the linear case, are discussed. For a special type of affine
Jacobi structures which are canonically exhibited (strongly-affine
or affine-homogeneous Jacobi structures) over a real vector space
of finite dimension, we describe the leaves of its characteristic
foliation as the orbits of an affine representation. These affine
Jacobi structures can be viewed as an analog of the
Kostant-Arnold-Liouville linear Poisson structure on the dual
space of a real finite-dimensional Lie algebra.
\end{abstract}
\begin{quote}
{\it Mathematics Subject Classification} (2000): 53D17, 53D05,
81S10.

\vspace{-5pt}

 {\it Key words and phrases}: Vector and affine
bundles, Jacobi manifolds, Lie algebroids.
\end{quote}

\vspace{-5pt}


\eject
\setcounter{section}{0}

\section{Introduction}
The Lie algebroid structures draw more and more attention in the
literature as structures generalizing the standard Cartan
differential calculus on differentiable manifolds.

There is a well-known correspondence between linear Poisson
structures on a vector bundle $\zp\colon A\rightarrow M$ and Lie
algebroid structures on the dual vector bundle $A^{\textstyle *}$,
which shows that the theory of Lie algebroids is, in fact, the
theory of linear Poisson brackets. This correspondence is built on
the fact that sections of $A^{\textstyle *} $ can be considered
linear functions on $A$. It can be easily extended to the
correspondence between {\it affine Poisson} structures on $A$
(whose brackets close on affine functions) and central extensions
of Lie algebroids, i.e., Lie algebroid structures on
$A^{\textstyle *} \times \R$ with the property that the section
$(0,1)$ is a central element of the Lie algebroid bracket. More
generally, a Lie algebroid structure on $A^{\textstyle *} \times
\R$  is uniquely represented by an {\it affine Jacobi} structure
on $A$ (whose Jacobi brackets close on affine functions) as well
as by a linear Poisson structure on $A\times \R$. Let us remark
that, as we will show later, the parts $\zL$ and $E$ of an affine
Jacobi tensor $(\zL,E)$ need not be affine themselves.

The presence of affine Poisson and Jacobi structures as
counterparts for Lie algebroids justifies  reconsideration of the
mentioned  relations in an affine setting, i.e., by admitting only
affine bundle structure on  $A$. The dual bundle is the vector
bundle $A^+=Aff(A,\R)$ whose fiber over $p\in M$ consists of
affine functions on the fiber $A_p$.  It has a distinguished
section corresponding to the constant function $1$ on $A$. In this
paper we will prove that there is one-to-one correspondence
between affine Jacobi structures on $A$ (the Jacobi bracket of
affine functions is an affine function) and Lie algebroids on
$A^+$.

The standard definition of a Lie algebroid structure on a vector
bundle $A$ consists of a Lie bracket defined on sections  and an
anchor map $\zr\colon A\rightarrow TM$. It is instructive to look
at a Lie algebroid as a restriction to sections of the
corresponding Schouten bracket $[\cdot,\cdot]_{SN}$ (which is, in
fact, a graded Poisson bracket) on the graded algebra of
multisections of $A$. The Schouten bracket  is graded
anticommutative, satisfies the graded Jacobi identity  and the
graded Leibniz rule. One can interpret a Poisson structure on $M$
as a canonical structure for the Schouten bracket of multivector
fields on $M$, i.e., as an element $\zL$ of Lie degree $-1$
satisfying the {\it master equation } $[\zL,\zL]_{SN}=0$. The
analogy with the classical Yang-Baxter equation is not an
accident, but the essence of the theory.

This point of view provides an easy passage to  the  theory of
Jacobi structures and Jacobi algebroids \cite{GM} (or generalized
Lie algebroids in the sense of \cite{IM2}). It is enough to
replace the (graded) Leibniz rule by the generalized Leibniz rule
which is valid for the first-order differential operators. The
obtained  bracket  is called  a  {\it Schouten-Jacobi bracket} on
$A$. Its restriction to sections of $A$ defines a Lie algebroid
structure on $A$, but its restriction to sections and functions
defines a Jacobi algebroid (\cite{GM, GM1}) structure on $A$ (a
generalized Lie algebroid in \cite{IM2}).

A Jacobi structure on $M$ in this setting turns out to be a {\it
canonical structure} (in the sense we explain later on)  for the
Schouten-Jacobi bracket of the  first order polydifferential
operators on $M$, i.e., skew-symmetric multidifferential
operators.

The passage from Poisson structures and Lie algebroids to Jacobi
structures and Jacobi algebroids is, essentially, the passage from
derivations to first-order differential operators. The notion of a
derivative depends on a reference frame (trivialization), but the
notion of  a first-order operator is not, and we make use of this
difference. A similar situation we encounter in physics when we
pass from a reference frame-dependent to frame-independent
description of a physical system. This is why frame-independent
formulations require affine bundles and Jacobi structures
(algebroids). We refer here to \cite{GGU, GGU1, MMS,SMM,MPL}
(time-dependent mechanics) and \cite{U}.

        With a Lie algebroid structure on $A$ we associate the
complete lift of multi-sections of $A$ to multivector fields on
$A$. It is a homomorphism of the Lie algebroid Schouten bracket
into the standard Schouten bracket and the complete lift of a
canonical structure of $A$ is a linear Poisson structure on $A$.

        Similarly, there is a complete lift  (cf. \cite {GM})
of a canonical structure for a Jacobi algebroid  $A$ to an affine
Jacobi structure  on $A$.

These remarks show that there is a need to look closer at affine
Jacobi brackets on affine (and also linear) bundles as to those
which are responsible for all these structures. This time,
however, the structures have an affine flavor.

The aim of this paper is a study of affine Jacobi structures on
affine and vector bundles and the corresponding Lie algebroids.

        A linear (resp., affine) Poisson structure on a vector bundle $A$
can be characterized by its behaviour with respect to the graded
algebra of polynomial functions on $A$ or with respect to the
Liouville (called also Euler) vector field  $\zD_A$ on $A$. Recall
that the Liouville vector field $\zD_A$ is the generator of the
one-parameter group of (positive) homoteties on $A$. For example,
a Poisson structure is linear, i.e., linear functions are closed
with respect to the Poisson bracket, if and only if one of the
following sentences is satisfied
\begin{itemize}
        \item the corresponding tensor $\zL$ is homogeneous with
respect to the Liouville vector field ($ \Ll_{\zD_A}\zL=-\zL$,
where ${\Ll}$ is the Lie derivative operator on $A)$.
        \item the Hamiltonian vector field of a linear function is linear.
\end{itemize}
A Poisson  structure is affine if and only $[Y,[X,\zL]]=0$ for
each pair of invariant vector fields $X, Y$ on $A$ (i.e., vertical
lifts of sections of $A$). In this case we say that $\zL$ is
affine homogeneous. This definition has its advantage, when
comparing with the action of the Liouville vector field, that it
can be used in non-commutative cases, i.e., for structures on Lie
groups or Lie groupoids.

        These characterizations cannot be extended to the case of
Jacobi structures. In particular, a linear Jacobi structure may
not be homogeneous and an affine Jacobi structure may not be
affine homogeneous, so one has to find the proper notion of
homogeneity in the affine case. In the paper we propose such
notion and we establish relations between different concepts and
specify the corresponding  Lie algebroids.

In Section~2 some definitions and results about Jacobi structures,
homogeneous multivectors in a vector bundle and Lie algebroids are
recalled. In Section~3 we discuss affine Jacobi structures on an
affine bundle in relation to linear Poisson structures on its
vector hull and Lie algebroids on the vector dual bundle. In
Section~4 we provide several examples. The most important ones are
given by the canonical structures which induce triangular
bialgebroid structures (Lie and Jacobi). In Section~5 we analyze
the relation between homogeneous and linear Jacobi structures
(previously, in Section~4, some results have been obtained).
Moreover, we introduce the notion of affine-homogeneous Jacobi
structures. We establish in Proposition~5.3 its relation to affine
and strongly-affine Jacobi structures. We remark that an affine
Jacobi structure is said to be strongly-affine if the hamiltonian
vector fields of affine functions are affine. On the other hand,
we prove that affine-homogeneous Jacobi structures on an affine
bundle $A$ correspond to Lie algebroids on the vector dual $A^+$
which have an ideal of sections of the subbundle spanned by $1_A$.
The Section~6 is devoted to the description of  leaves of the
characteristic foliation of a strongly-affine  Jacobi structure on
a vector space, as the orbits of an affine representation of a Lie
group on the vector space. It can be viewed as a generalization of
viewing symplectic leaves of the Kostant-Arnold- Liouville linear
Poisson structure on the dual space to a Lie algebra as the orbits
of the coadjoint action of the corresponding Lie group.
\setcounter{equation}{0}

\section{Jacobi manifolds and Lie algebroids}\label{seccion2}

A {\it Jacobi manifold} \cite{L2} is a differentiable manifold $M$
endowed with a pair ({\it Jacobi structure})  $(\Lambda,E)$, where
$\Lambda$ is a $2$-vector and $E$ is a vector field on $M$
satisfying
\[
[\Lambda,\Lambda]_{SN}=-2E\wedge \Lambda,\makebox[1cm]{}
[E,\Lambda]_{SN}=0.
\]
Here $[\cdot,\cdot]_{SN}$ denotes the Schouten bracket. Note that
we use the version of  the  Schouten-Nijenhuis  bracket  which
gives a graded  Lie  algebra  structure  on  multivector  fields
and which differs from the classical one (\cite{BV,V1}) by signs.
For this type of manifolds, a bracket of functions (the {\it
Jacobi bracket}) is defined by
\[
\{f,g\}_{(\Lambda,E)}=\Lambda(df,dg)+fE(g)-gE(f),
\] for all
$f,g\in C^\infty(M,\R).$ This bracket is skew-symmetric, satisfies
the Jacobi identity and it is a first-order differential operator
on each of its arguments, with respect to the ordinary
multiplication of functions.  We will  often identify  the Jacobi
bracket with the first-order bidifferential operator $\zL+I\we E$,
where $I$ is the identity on $C^\infty(M,\R)$. The space
$C^\infty(M,\R)$ of $C^\infty$ real valued functions on $M$
endowed with the Jacobi bracket is a local Lie algebra (see
\cite{Ki}). Conversely, a local Lie algebra on $C^\infty(M,\R)$
defines a Jacobi structure on $M$ (see \cite{GL,Ki}). Note that
Poisson manifolds \cite{L1} are Jacobi manifolds  with $E=0$.

Other interesting examples of Jacobi manifolds are contact and
locally conformal symplectic manifolds (see for example
\cite{GL,L2,Vai1} for the definition of these types of manifolds).


For a Jacobi manifold $(M,\Lambda,E)$, one can consider the
homomorphism of $C^\infty(M,\R)$-modules
$\#_{\Lambda}:\Omega^1(M)\to {\mathfrak X}(M)$ given by
$\#_{\Lambda}(\alpha)(\beta)=\Lambda(\alpha,\beta),$ for all
$\alpha,\beta\in\Omega^1(M).$ If $f$ is a $C^\infty$ real-valued
function on a Jacobi manifold $M$, then the vector field
$X_f^{(\Lambda,E)}$ defined by $ X^{(\Lambda,E)}_f=\#_\Lambda(df)
+ fE $ is called {\it the hamiltonian vector field } associated
with $f.$ It should be noticed that the hamiltonian vector field
associated with the constant function $1$ is just $E$.

Now, for every $x\in M,$ we consider the subspace ${\cal
F}^{(\Lambda,E)}_x$ of $T_xM$ generated by all the hamiltonian
vector fields evaluated at the point $x.$ In other words, ${\cal
F}^{(\Lambda,E)}_x=(\#_\Lambda)_x(T_x^*M) + <E_x>.$ Since ${\cal
F}^{(\Lambda,E)}$ is involutive  and  finitely generated, it is
easily seen that ${\cal F}^{(\Lambda,E)}$ defines a generalized
foliation in the sense of Sussmann \cite{Sus}, which is called
{\it the characteristic foliation } (see \cite{DLM,GL}). Moreover,
the Jacobi structure on $M$ induces a Jacobi structure on each
leaf. In fact, if $L$ is the leaf over a point $x$ of $M$ and if
$E_x\notin Im(\#_\Lambda)_x$ (or equivalently, the dimension of
$L$ is odd) then $L$ is a contact manifold. If $E_x\in
Im(\#_\Lambda)_x$ (or equivalently, the dimension of $L$ is even),
then L is a l.c.s. manifold (see \cite{GL,Ki}; see also
\cite{IM3}).
\begin{definition}
A tensor field $X$ on $M$ is {\it homogeneous of degree $k$} with
respect to a vector field $\zD$ if $\Ll_\zD X=kX$, where $\Ll$ is
the Lie derivative. If $k=1-n$ for a contravariant $n$-tensor
field $X$, then we will call $X$ just {\it homogeneous}. In
particular, a bivector field is homogeneous if $\Ll_\zD X=-X$.
\end{definition}

We will often identify sections $\zm$ of the dual bundle $A^*$
with linear (along fibres)  functions  $\zi_\zm$  on the vector
bundle $A$: $\zi_\zm(X_p)=<\zm(p),X_p>$. Note that if $f:A\to \R$
is a smooth real function and $\Delta_A$ is the Liouville vector
field of $A$ then
\begin{equation}\label{lineal}
f \mbox{ is linear} \Leftrightarrow \Delta_A(f)=f.
\end{equation}
We recall that $\Delta _A$ is the vector field on $A$ given by
$\Delta _A=\displaystyle \sum_{\alpha}y^\alpha
\frac{\partial}{\partial y^\alpha }$, for fibred coordinates
$(x^i, y^\alpha)$.

 On the other hand, a $2$-vector $\zL$ on $A$ is {\it linear
} if and only if the induced bracket is closed on linear
functions, that
is,$<\zL,\D\zi_{\zm}\we\D\zi_{\zn}>=\{\zi_\zm,\zi_\zn\}_\zL $ is
again a linear function associated with  an  element
$[\zm,\zn]_\zL$. The  operation $[\zm,\zn]_\zL$ on sections  of
$A^*$ is called the {\it bracket induced by $\zL$}. If $\zL$ is a
linear $2$-vector field on $A$ and $f,g:A\to \R$ are basic
functions then
\begin{equation}\label{basica}
\begin{array}{l}
<\zL,\D\zi_{\zm}\we\D f> \mbox{ is a basic function
}\makebox[2cm]{and}<\zL, \D f\wedge\D g>=0.
\end{array}
\end{equation}
Using the above facts, it is easy to prove that $\zL$ is linear if
and only if it is a homogeneous bivector field on $A$ with respect
to $\Delta_A$.

\medskip
A {\it Lie algebroid structure} \cite{Mk} on a differentiable
vector bundle $\pi:A\to M$ is a pair which consists of a Lie
algebra structure $\lcf \cdot,\cdot\rcf$ on the space $\Gamma(A)$
of the global sections of $\pi:A\to M$ and a homomorphism of
vector bundles $\rho:A\to TM$, the {\it anchor map}, such that if
$\rho:\Gamma(A)\to {\mathfrak X}(M)$ also denotes the homomorphism
of $C^\infty(M,\R)$-modules induced by the anchor map, then
\[
\lcf X,fY\rcf=f\lcf X,Y\rcf + \rho(X)(f)Y,
\]
for all $X,Y\in \Gamma(A)$ and $f\in C^\infty(M,\R)$. It follows
that $\rho:(\Gamma(A),\lcf\cdot,\cdot\rcf)\to ({\mathfrak
X}(M),[\cdot,\cdot])$ is a Lie algebra homomorphism. Note that the
anchor is uniquely determined by the Lie algebroid bracket.

\begin{theorem}\cite{Co} There is a one-one correspondence between  Lie
algebroid brackets $\lna\cdot,\cdot\rna_\zL$ on the vector bundle
$A$ and  homogeneous (linear) Poisson structures $\zL$ on the dual
bundle $A^*$ determined by
\[
\zi_{\lna X,Y\rna_\zL}=\{\zi_X,\zi_Y\}_\zL=\zL(d\zi_X,d\zi_Y).
\]
\end{theorem}
Every  Poisson  structure  $\zL$  on  $M$  determines  a  Lie
algebroid bracket $\lna\cdot,\cdot\rna_\zL$ on $T^*M$ with the
anchor $\#_\zL$ and the bracket $\lna\cdot,\cdot\rna_\zL$ defined
by $\lna\za,\zb\rna_\zL=i_{\#_\Lambda(\alpha)}d\beta-
i_{\#_\Lambda(\beta)}d\alpha +d(\Lambda(\alpha,\beta)).$ Also
every Jacobi manifold $(M,\Lambda,E)$ is associated with a Lie
algebroid $(T^*M\times \R, \lcf\cdot,\cdot\rcf_{(\Lambda,E)},$ $
\tilde{\#}_{(\Lambda,E)}),$ where the Lie bracket
$\lcf\cdot,\cdot\rcf_{(\Lambda,E)}:(\Omega^1(M)\times
C^\infty(M,\R))^2\to \Omega^1(M)\times C^\infty(M,\R)$, and the
anchor map $\widetilde{\#}_{(\Lambda,E)}:\Omega^1(M)\times
C^\infty(M,\R)\to {\mathfrak X}(M)$ are defined by (see \cite{KS})
\begin{equation}\label{e0}
\begin{array}{rcl}
\kern-10pt\lcf(\alpha,f),(\beta,g)\rcf_{(\Lambda,E)}\kern-10pt&=&
\kern-10pt(i_{\#_\Lambda(\alpha)}d\beta-i_{\#_\Lambda(\beta)}d\alpha
+d(\Lambda(\alpha,\beta)) +
f\Ll_{E}\beta-g\Ll_{E}\alpha-i_E(\alpha\wedge \beta),\\&&
\Lambda(\beta,\alpha) + \#_\Lambda(\alpha)(g)
-\#_\Lambda(\beta)(f) + fE(g)-gE(f)),\\
\widetilde{\#}_{(\Lambda,E)}(\alpha,f)&=&\#_\Lambda(\alpha) + fE,
\end{array}
\end{equation}
for $(\alpha,f),(\beta,g)\in \Omega^1(M)\times C^\infty(M,\R).$


As it has been observed in \cite{KoS}, a {\it Lie algebroid}
structure on a vector bundle $A$ can be identified  with  a  {\it
Gerstenhaber algebra} structure (in the terminology of \cite{KoS})
on the exterior algebra of multisections of $A$, $\Gamma(\wedge
A)=\oplus_{k\in \Z}\Gamma(\wedge^kA),$ which is just a {\it graded
Poisson bracket } ({\it Schouten bracket}) on $\Gamma(\wedge A)$
of degree -1 (linear).

A Schouten bracket induces  the  well-known generalization of the
standard Cartan calculus of differential forms and vector fields
\cite{Mk}. The exterior derivative $\D:\Gamma(\wedge^k A)\to
\Gamma(\wedge^{k+1}A)$ is defined by the standard formula \bea\nn
        \D\zm(X_1,\dots,X_{k+1}) &=& \sum_i (-1)^{i+1}
\lna X_i,\zm(X_1,\dots,\widehat{X}_i,\dots ,X_{k+1})\rna \\ &+&
\sum_{i<j} (-1)^{i+j}\zm (\lna X_i,X_j\rna, X_1,\dots ,
\widehat{X}_i, \dots ,\widehat{X}_j, \dots , X_{k+1}).\label{D}
\eea

For $X \in\zG(A)$, the contraction $i_X \colon \Gamma(\wedge^p A)
\rightarrow\Gamma(\wedge ^{p-1}A)$ is defined in the standard way
and the Lie differential operator $ \Ll_X$
         is defined by the graded commutator
$\Ll_X = i_X \circ \D +\D \circ i_X.$

\medskip
It is obvious  that  the   notion of Schouten bracket  extends
naturally to more general gradings in the algebra. For  a  graded
commutative algebra   with   unity   $\1$,   a   natural
generalization of  a graded  Poisson  bracket  is  {\it  graded
Jacobi bracket}. The only difference is that we replace the
Leibniz rule by the generalized Leibniz rule: \be\label{first}
\lcf a,bc\rcf=\lcf  a,b\rcf c+(-1)^{(\vert a\vert+k)\vert
b\vert}b\lcf a,c\rcf-\lcf a,\1\rcf bc. \ee Graded Jacobi brackets
on $\Gamma(\wedge A)$  of degree $k=-1$ (linear) is called  {\it
Schouten-Jacobi} brackets. An element $X\in\Gamma(\wedge^2 A)$  is
called  a  {\it canonical structure} for a Schouten or
Schouten-Jacobi bracket $\lna\cdot,\cdot\rna$ if $\lna X,X\rna=0$.

Since Schouten brackets  on  $\Gamma(\wedge A)$  are just Lie
algebroid structures on $A$ (see \cite{KoS}), by  a {\it
generalized Lie algebroid} (or {\it Jacobi algebroid}) structure
on $A$ we mean a Schouten-Jacobi bracket on $\Gamma(\wedge A)$.
The generalized Lie algebroids are in one-one correspondence with
pairs consisting of a Lie algebroid $A$ and a $1$-cocycle
$\phi_0\in \Gamma(A^*)$ relative to the Lie algebroid exterior
derivate $d$, i.e., $d\phi_0=0$, (cf. \cite{GM,GM1,IM2}).

A canonical example  of  a Jacobi algebroid is $(\zT M, (0,1))$,
where $\zT M=TM\oplus\R$ is the  Lie algebroid of first-order
linear differential operators on $C^\infty(M,\R)$ with the bracket
\be\label{ddd} [(X,f),(Y,g)]_1=([X,Y],X(g)-Y(f)),\quad X,Y\in
\X(M),\quad f,g\in C^\infty(M,\R), \ee and  the  1-cocycle
${\phi_0}=(0,1)$ is ${\phi_0}((X,f))=f$. Note  that  we  have  the
canonical decomposition $X=X_1+I\we X_2$ of any tensor
$X\in\zG(\wedge^k\zT M)$, where $X_1$ (resp. $X_2$) is a
$k$-vector  field  (resp. $(k-1)$-vector field) and  $I$
represents  the  identity operator on  $C^\infty(M,\R)$ which is a
generating section of $\R$ in $TM\oplus\R$.
   A   canonical
structure with respect to the corresponding  Schouten-Jacobi
bracket  on the Grassmann algebra  $\zG(\wedge \zT M)$  of
first-order polydifferential operators on $C^\infty(M,\R)$, which
we will denote by $[\cdot,\cdot]_1$, turns out to be a standard
Jacobi structure. Indeed, it is  easy  to  see that the
Schouten-Jacobi bracket reads \[\begin{array}{rcl} [A_1+I\we
A_2,B_1+I\we B_2]_1 &=&[A_1,B_1]_{SN} +(-1)^aI\we
[A_1,B_2]_{SN}+I\we[A_2,B_1]_{SN}\\ &&+ aA_1\we B_2-(-1)^abA_2\we
B_1+(a-b)I\we A_2\we B_2.\end{array}\] Hence, the bracket
$\{\cdot,\cdot\}$ on $C^\infty(M,\R)$  defined by a bilinear
differential operator $\zL+I\we\zG\in\zG(\wedge^2\zT
M)=\zG(\wedge^2(TM)\oplus\zG(TM))$ is  a Lie bracket (Jacobi
bracket on $C^\infty(M,\R)$) if and only if
\[
[\zL+I\we\zG,\zL+I\we\zG]_1=[\zL,\zL]_{SN}+2I\we
[\zG,\zL]_{SN}+2\zL\we\zG=0.
\]
Thus, we get the conditions defining a Jacobi structure on $M$.

There is another  approach  to  Lie  algebroids.  As it was shown
in \cite{GU1, GU2}, a Lie algebroid   structure (or the
corresponding Schouten  bracket)  is  determined  by the Lie
algebroid   lift $X\mapsto X^c$ which  associates  with
$X\in\zG(\wedge A)$ a multivector field $X^c$ on $A$. The complete
lifts of  Lie  algebroids  are  described  as follows. For a given
Lie  algebroid  structure  on a vector bundle $A$ over $M$ there
is a unique {\it complete lift} of elements
$X\in\zG(\bigwedge^kA)$ of the Grassmann algebra $\Gamma(\wedge
A)=\oplus_k \zG(\bigwedge^kA)$  to multivector fields
$X^c\in\zG(\bigwedge^k(TA))$ on $A$ , such that
\[f^c=\zi_{\D f}, \;\;\; X^c(\zi_\zm)=\zi_{\Ll_X\zm},\;\;\;
(X\we Y)^c=X^c\we Y^v+X^v\we Y^c, \] for $f\in C^\infty(M,\R),$
$X,Y\in \zG(A)$ and $  \zm\in \zG(A^*)$, where $X\mapsto  X^v$  is
the standard vertical lift of tensors  from  $\Gamma(\wedge A)$ to
tensors  from $\Gamma(\wedge TA)$. Moreover,  this  complete  lift
is a   homomorphism of   the corresponding Schouten brackets: \[
\lna X,Y\rna^c=[X^c,Y^c]_{SN}
 \makebox[2cm]{and}
\lna X,Y\rna^v=[X^c,Y^v]_{SN}. \] Note that in the particular case
$A=TM$, the above complete lift reduces to well-known tangent lift
of multivector fields on  $M$ to $TM$ (cf. \cite{GU,GU1,IY,MX1}).

Recently, in \cite{Gr} the notion of {\it a Lie QD-algebroid} has
been introduced as a  vector bundle $A$ on a manifold $M$ endowed
with  a Lie algebra bracket $[\cdot,\cdot]$  on the
$C^\infty(M,\R)$-module of sections of  $A$ and a
$C^\infty(M,\R)$-linear map $X\in \Gamma(A)\to \hat{X}\in {\frak
X}(M)$  which satisfy $[X,fY]=f[X,Y] + \hat{X}(f)Y$ for all
$X,Y\in \Gamma(A)$, $f\in C^\infty(M,\R)$. If $A$ is a line bundle
then a Lie QD-algebroid is a local Lie algebra in the sense of
Kirillov \cite{Ki}. If $rank\; A>1,$ a Lie QD-algebroid on $M$ is
just a Lie algebroid on $M$ (for more details, see \cite{Gr}).


\section{Affine Jacobi structures on affine bundles}
\setcounter{equation}{0} Let $\pi:A\to M$ be an affine bundle over
$M$ of rank $n$ modelled on a vector bundle $V(\pi):V(A)\to M$,
that is, $\pi:A\to M$ is a (locally trivial) smooth bundle such
that the fiber at the point $x\in M$, $A_x=\pi^{-1}(x),$ is an
affine space modelled on the vector space $V_x=(V(\pi ))^{-1}(x)$,
and we pass from one local trivialization to another using the
group of affine transformations. If $p\in M,$ we denote by
$Aff(A_p,\R)$ the vector space of affine functions from the fiber
$A_p$  of $\pi:A\to M$ at $p$ and by $A^+$ the vector bundle
$A^+=\bigcup_{p\in M}Aff(A_p,\R)\to M$ of rank $n+1$. Note that
$A^+$ has a distinguished $1$-section, $\tilde{1}:M\to A^+$,
defined by the constant function $1$ on $A$, i.e.,
$\tilde{1}(p)=1_{{A_p}}\in Aff(A_p,\R).$ We call
$A^\dag=(A^+,\tilde 1)$ {\it the (special) vector dual}  of  the
affine bundle $A$. In general, by a {\it special vector bundle}
(cf. \cite{GGU,GGU1}) we  mean a vector bundle with a
distinguished nowhere-vanishing section, so that $A\mapsto A^\dag$
gives rise to a  (contravariant)  functor  from the category of
affine bundles to the category of special vector bundles.

We also have a dual functor which assigns to every special vector
bundle $(V,X)$ an affine  bundle  $(V,X)^\ddag$  which  is the
affine  bundle defined as the $1$-level set of the linear function
$\zi_X$ in the dual vector bundle $V^*$. It is  easy  to see that
for an affine bundle $A$ we have $(A^\dag)^\ddag\simeq A$,  so
that we can identify $A$ with an affine subbundle of $\hat
A=(A^+)^*$ (in fact, $A=\iota_{\tilde{1}}^{-1}(1)).$ Note that we
have a full duality, since also $(V,X)=(((V,X)^\ddag)^+,\tilde
1)$.

Using this fact, one can prove that there is a one-to-one
correspondence between affine functions on $A$ and linear
functions on $\hat{A}$. In fact, if $a:A\to \R$ is an affine
function on each fiber of $A$ then the corresponding linear
function $\bar{a}:\hat{A}\to \R$ on each fiber of  $\hat{A}$ is
given by
$
\bar{a}(\psi_p)=\psi_p(a_{|A_p}),$ for all $\psi_p\in
\hat{A}_p=(A_p^+)^\ast. $ Note that $\bar{a}_{|A}=a.$

Moreover, there is an obvious natural one-to-one correspondence
between affine functions and sections of $A^+$ which associates
with the section of $A^+$, $\tilde{a}\in \Gamma(A^+)$, the affine
function $a:A\to \R$ and the linear map
$\zi_{\tilde{a}}:\hat{A}\to \R$ is just the function $\bar{a}.$

\begin{definition}
A Jacobi structure on a vector bundle (resp. affine bundle) is
called {\it linear} (resp. {\it affine}) if  the corresponding
Jacobi bracket  of linear functions is again a linear function
(resp. the bracket of affine functions is an affine function).
\end{definition}

Now, we consider   an affine Jacobi structure
$(\Lambda_A,E_A)$ on an affine bundle $A.$ Denote by
$\{\cdot,\cdot\}_{(\Lambda_A,E_A)}$ the corresponding Jacobi
bracket.

If the rank of $A$ is zero, i.e., $A=M\times \{x_0\}$, then
$(\Lambda_A,E_A)$ induces a Jacobi structure $(\Lambda,E)$ over
$M$.

Moreover,  $A^+=M\times \R$ and thus $\Gamma(A^+)\cong
C^\infty(M,\R).$ Therefore, the Jacobi bracket
$\{\cdot,\cdot\}_{(\Lambda,E)}$ induces a Lie QD-algebroid
structure on $A^+$. Conversely, if $\lcf\cdot,\cdot\rcf$ defines a
Lie QD-algebroid structure on $A^+$ then we have that a local Lie
algebra structure on the real line bundle $A^+=M\times \R\to M$ or
equivalently, a Jacobi structure on $M$, i.e., an affine Jacobi
structure on $A.$

Now, we suppose that the rank of the affine bundle $A$ is
non-zero. Then, we have the following result.
\begin{lemma}\label{l1}
Let $f:A\to \R$ be a basic function.
\begin{enumerate}
\item If $a:A\to \R$ is an affine function, then
$\{f,a\}_{(\Lambda_A,E_A)}-f\{1,a\}_{(\Lambda_A,E_A)}$  is a basic
function.
\item If $g:A\to \R$ is a basic function, then
$\{f,g\}_{(\Lambda_A,E_A)}$ is a basic function and
\[
\{f,g\}_{(\Lambda_A,E_A)}=f\{1,g\}_{(\Lambda_A,E_A)}+g\{f,1\}_{(\Lambda_A,E_A)
}.
\]
\end{enumerate}
\end{lemma}
{\bf Proof.-} Let $p$ be a point of $M$. Since $rank \; A>0$ one
can choose an affine function  $b:A\to \R$ such that the linear
function associated with the affine function $b_{|A_p}:A_p\to \R$
is non-zero. Then,
\begin{equation}\label{r1}
\{bf,a\}_{(\Lambda_A,E_A)}=b\{f,a\}_{(\Lambda_A,E_A)} +
f\{b,a\}_{(\Lambda_A,E_A)} - bf\{1,a\}_{(\Lambda_A,E_A)}.
\end{equation}
Since $a$ is an affine function, then $\{bf,a\}_{(\Lambda_A,E_A)}$
and $f\{b,a\}_{(\Lambda_A,E_A)}$ are affine functions and
therefore, from (\ref{r1}), we have that
$b(\{f,a\}_{(\Lambda_A,E_A)}-f\{1,a\}_{(\Lambda_A,E_A)})$ is
affine, that is,
$(\{f,a\}_{(\Lambda_A,E_A)}-f\{1,a\}_{(\Lambda_A,E_A)})_{|A_p}$ is
a constant function. Therefore,
$\{f,a\}_{(\Lambda_A,E_A)}-f\{1,a\}_{(\Lambda_A,E_A)}$ is a basic
function, i.e., {\it (i)} holds.

\medskip
If $g$ is a basic function then, $\{bf,g\}_{(\Lambda_A,E_A)}$ and
$f\{b,g\}_{(\Lambda_A,E_A)}$ are affine functions. Moreover, using
{\it (i)} for the affine function $a\equiv 1$ and the basic
function $g$, we have that $\{1,g\}_{(\Lambda_A,E_A)}$ is a basic
function. Thus, since
\[
\{bf,g\}_{(\Lambda_A,E_A)}=b\{f,g\}_{(\Lambda_A,E_A)} +
f\{b,g\}_{(\Lambda_A,E_A)} - bf\{1,g\}_{(\Lambda_A,E_A)},
\]
the function $b\{f,g\}_{(\Lambda_A,E_A)}$ is affine. Consequently,
$(\{f,g\}_{(\Lambda_A,E_A)})_{|A_p}$ is a constant function. This
proves that $\{f,g\}_{(\Lambda_A,E_A)}$ is a  basic function.
Furthermore, from $(i),$ we obtain that
\[
\begin{array}{rcl}
\{g,bf\}_{(\Lambda_A,E_A)}-g\{1,bf\}_{(\Lambda_A,E_A)}&=&f(\{g,b\}_{(\Lambda_A
,E_A)} -g\{1,b\}_{(\Lambda_A,E_A)})
\\&&-b(\{f,g\}_{(\Lambda_A,E_A)}-f\{1,g\}_{(\Lambda_A,E_A)}
-g\{f,1\}_{(\Lambda_A,E_A)})\end{array}\]  is a basic function,
which implies that (see $(i)$)
$b(\{f,g\}_{(\Lambda_A,E_A)}-f\{1,g\}_{(\Lambda_A,E_A)}-g\{f,1\}_{(\Lambda_A,E
_A)})$ is a basic function, that is,
$(\{f,g\}_{(\Lambda_A,E_A)}-f\{1,g\}_{(\Lambda_A,E_A)}-g\{f,1\}_{(\Lambda_A,E_A)
})_{|A_p}=0.$ Therefore,
$\{f,g\}_{(\Lambda_A,E_A)}=f\{1,g\}_{(\Lambda_A,E_A)}+g\{f,1\}_{(\Lambda_A,E_A)
}.$\hfill$\Box$

Now, we will describe  a Lie algebroid structure
$(\lcf\cdot,\cdot\rcf^{+},\rho^{+})$ on the vector bundle $A^+$
using the natural one-to-one correspondence between sections of
$A^+$ and the space of the affine functions on $A.$ In fact, if we
denote by $\tilde{a}$ the section of $A^+$ associated with the
affine function $a:A\to \R$, then, we define the pair
$(\lcf\cdot,\cdot\rcf^{+},\rho^{+})$ on $A^+$ as follows
\begin{equation}\label{r2}
\left.\begin{array}{rcl}
\lcf\tilde{a},\tilde{b}\rcf^{+}&=&\widetilde{\{a,b\}}_{(\Lambda_A,E_A)},
\\[5pt]
\rho^{+}(\tilde{a})({f_M})\circ \pi&=&\{a,{f_M}\circ
\pi\}_{(\Lambda_A,E_A)}-({f_M}\circ \pi)\{a,1\}_{(\Lambda_A,E_A)},
\end{array}\right\}
\end{equation}
for all $a,b:A\to \R$ affine functions and ${f}_M\in
C^\infty(M,\R).$
\begin{theorem}\label{p1}
Let $(\Lambda_A,E_A)$ be an affine Jacobi structure on an affine
bundle $\pi:A\to M$ and assume that the rank of $A$ is $n$, $n>0$.
Then, the bracket $\lcf\cdot,\cdot\rcf^{+}:\Gamma(A^+)\times
\Gamma(A^+)\to \Gamma(A^+)$ and the map $\rho^{+}:\Gamma(A^+)\to
{\mathfrak X}(M)$ given as in (\ref{r2}) define a Lie algebroid
structure on $A^+$.
\end{theorem}
{\bf Proof.-} Since $\{\cdot,\cdot\}_{(\Lambda_A,E_A)}$ is
skew-symmetric and it satisfies the Jacobi identity, one deduces
easily that $(\Gamma(A^+),\lcf\cdot,\cdot\rcf^{+})$ is a Lie
algebra. Moreover, using Lemma \ref{l1} and the fact that
$\{\cdot,\cdot\}_{(\Lambda_A,E_A)}$ is a first-order
bi-differential operator, we obtain that $\rho^{+}:A^+\to TM$ is a
homomorphism of vector bundles.

Finally, from the fact that $\{\cdot,\cdot\}_{(\Lambda_A,E_A)}$ is
a first-order bidifferential operator, we conclude that
$(\lcf\cdot,\cdot\rcf^{+},$ $\rho^{+})$ is a Lie algebroid
structure on $A^+$. \hfill$\Box$

\begin{remark}\label{remark3.2'}{\rm Let $(\Lambda_A,E_A)$ be an affine Jacobi
structure on an affine bundle $\pi:A\to M$ of rank $n,$ $n>0$, and
$(\lcf\cdot,\cdot \rcf^+,\rho^+)$ be the corresponding Lie
algebroid structure on $A^+$. Denote by $\bar\Lambda_{\hat{A}}$
the linear Poisson structure on $\hat{A}=(A^+)^*$ induced by the
Lie algebroid structure $(\lcf\cdot,\cdot\rcf^+,\rho^+).$ Then, we
have that
\begin{equation}\label{3.2 0}
\iota_{\lcf\tilde{a},\tilde{b}\rcf^+}=
\{\bar{a},\bar{b}\}_{\bar\Lambda_{\hat{A}}},\end{equation} for
$a,b:A\to \R$ affine functions on $A$. On the other hand, if
$\Delta_{\hat{A}}$ is the Liouville vector field we deduce that $
(\Delta_{\hat{A}}(\iota_{\tilde{1}}))_{|A}=(\iota_{\tilde{1}})_{|A}=1,
$
which implies that $\Delta_{\hat{A}}$ is a transverse vector field
of $A$ as a submanifold of $\hat{A}$. Thus, using (\ref{lineal}),
Proposition 2.3 in \cite{DLM} and since $\bar\Lambda_{\hat{A}}$ is
a homogeneous Poisson structure on $\hat{A},$ we obtain that there
exists a Jacobi structure $(\Lambda'_A,E'_A)$ on $A$ such that
\begin{equation}\label{3.2 1}
\{a,b\}_{(\Lambda'_A,E'_A)}=(\{\bar{a},\bar{b}\}_{\bar\Lambda_{\hat
A}})_{|A},
\end{equation}
for $a,b:A\to \R$ affine functions on $A.$ Therefore, from
(\ref{r2}), (\ref{3.2 0}) and (\ref{3.2 1}), we conclude that
$\{a,b\}_{(\Lambda'_A,E'_A)}=\{a,b\}_{(\Lambda_A,E_A)},$ i.e.,
$(\Lambda'_A,E'_A)$ is just the affine Jacobi structure
$(\Lambda_A,E_A).$

}
\end{remark}

 Now, we will prove the
converse of Theorem \ref{p1}.
\begin{theorem}\label{p2}
Let $(\lcf \cdot,\cdot\rcf^{+},\rho^{+})$  be a Lie algebroid
structure on $\pi^+:A^+\to M$ and assume that the rank of $A$ is
$>0$. Then, there exists a unique affine Jacobi structure
$(\Lambda_A,E_A)$ on $A$ such that
\begin{equation}\label{r3}
\widetilde{\{a,b\}}_{(\Lambda_A,E_A)}=\lcf\tilde{a},\tilde{b}\rcf^{+},
\;\;\; \forall a,b:A\to \R\mbox{ affine functions}.
\end{equation}
\end{theorem}

{\bf Proof.-} The uniqueness is deduced from the fact that a
Jacobi structure is characterized by the Jacobi bracket of linear
functions and the Jacobi bracket of a linear function and the
constant function $1$. Thus, two Jacobi structures satisfying
(\ref{r3}) are equal.

\medskip

Now, we will define a Jacobi structure on $A$ which satisfies
(\ref{r3}). Denote by $\bar\Lambda_{\hat{A}}$ the linear Poisson
structure on $\hat{A}$ induced by the Lie algebroid structure
$(\lcf\cdot,\cdot\rcf^{+},\rho^{+})$. Then, proceeding as in
Remark \ref{remark3.2'}, we have that there exists a Jacobi
structure $(\Lambda_A,E_A)$ on $A$ such that
\begin{equation}\label{e3.4}
\{a,b\}_{(\Lambda_A,E_A)}=(\{\bar{a},\bar{b}\}_{\bar\Lambda_{\hat{A}}})_{|A},
\end{equation}
for $a,b:A\to \R$ affine functions on $A$. In fact, if
$\Delta_{\hat{A}}$ is the Liouville vector field on $\hat{A}$,
$E_{\hat{A}}$  is the hamiltonian vector field of the linear
function $\iota_{\tilde{1}}:\hat{A}\to \R$ with respect to
$\bar\Lambda_{\hat A}$ and $\Lambda_{\hat{A}}$ is the $2$-vector
on $\hat{A}$ given by
\begin{equation}\label{r5}
\Lambda_{\hat{A}}=\bar{\Lambda}_{\hat{A}}-\Delta_{\hat{A}}\wedge
E_{\hat{A}}, \end{equation}
 then,  we obtain
that $(\Lambda_{\hat{A}}, E_{\hat{A}})$ is a Jacobi structure on
$\hat{A}$ and, from (\ref{e3.4}), it follows  that $\Lambda_A$
(respectively, $E_{A}$) is the restriction to $A$ of the
$2$-vector $\Lambda_{\hat{A}}$ (respectively, $E_{\hat{A}}$).
Moreover, using again (\ref{e3.4}), it follows that the Jacobi
structure $(\Lambda_A,E_A)$ is affine and, in addition, (\ref{r3})
holds. \hfill$\Box$

\begin{remark}
{\rm Let $\pi : A\to M$ be an affine bundle such that  the dual
vector bundle $\pi ^+:A^+\to M$ carries a Lie algebroid structure
$(\lcf \cdot ,\cdot \rcf ^{+},\rho ^{+})$. Let $(x^l)_{l=1,\dots
,m}$ be local coordinates on an open subset $U_M$ of  $M$. We can
consider $\{e_0, {e}_1,\dots, {e}_n\}$ a local basis of sections
of $\pi^+:A^+\to M$ such that ${e}_0$ is the section of $A^+$
associated with the linear  function $\zi_{\tilde 1}:\hat{A}\to
\R$. Then, there is  an open coordinate neighbourhood $U$ on
$\hat{A}$ with coordinates $(x^1,\dots ,x^m,\zi_{\tilde
1},y^1,\dots ,y^n)$, where $y^\alpha:U\to \R$ is the linear map
associated with $e_\alpha$. If the structure functions of
$\lcf\cdot,\cdot\rcf^{+}$ and the components of the anchor map
$\rho^+$ for these coordinates are
$c_{\alpha\beta}^\gamma,\rho_\alpha^l\in C^\infty(U_M,\R)$ then
$(\Lambda _A,E_A)$ is given by
\begin{equation}\label{e5}\begin{array}{rcl}
\Lambda_A&=&\displaystyle\sum_{\mbox{ \tiny ${\begin{array}{c}{\alpha<\beta}\\
{ \alpha,\beta=1,\dots ,n}\end{array}}$}}
\kern-10pt\sum_{\gamma=1,\dots ,n} (c_{0\alpha}^\gamma y^\gamma
y^\beta -c_{0\beta}^\gamma y^\gamma y^\alpha +
c_{\alpha\beta}^\gamma y^\gamma+c_{0\alpha}^0 y^\beta -
c_{0\beta}^0 y^\alpha + c_{\alpha\beta}^0) \frac{\partial
}{\partial y^\alpha}\wedge \frac{\partial }{\partial y^\beta} \\&&
+ \displaystyle\sum_{l=1,\dots ,m}
\displaystyle\sum_{\alpha=1,\dots ,n} (\rho_{\alpha}^l -
y^\alpha\rho_{0}^l)\frac{\partial }{\partial y^\alpha}\wedge
\frac{\partial }{\partial x^l}, \\ E_A
&=&\displaystyle\sum_{\beta=1,\dots ,n}(
\displaystyle\sum_{\gamma=1,\dots ,n}c_{0\beta}^\gamma y^\gamma +
c_{0\beta}^0)\frac{\partial }{\partial y^\beta} +
\displaystyle\sum_{l=1,\dots ,m} \rho_{0}^l \frac{\partial
}{\partial x^l}. \end{array}\end{equation}

Note that,  in general, the local components of $\Lambda_A$ are
not affine functions. }
\end{remark}
From Theorems \ref{p1} and \ref{p2} and taking into account that
there exists a one-to-one correspondence between affine Jacobi
structures on an affine bundle $A$ of rank zero and Lie
QD-algebroid structures on $A^+,$  we obtain that

\begin{corollary}\label{c1}
Let $\pi:A\to M$ be an affine bundle of rank $n$. Then:
\begin{enumerate}
\item If $n>0,$ there is a one-to-one correspondence between affine Jacobi
brackets on  $\pi:A\to  M$ and Lie algebroid structures on the
vector bundle $A^+$ uniquely determined by the equation
(\ref{r3}).
\item If $n=0$, there is a one-to one correspondence between affine Jacobi
brackets on  $\pi:A\to  M$ and local Lie algebra structures  on
$A^+=M\times \R.$
\end{enumerate}
\end{corollary}

Using the equivalence between Lie algebroids and linear Poisson
brackets we can formulate also a Poisson version of the above
Corollary. First, we introduce the following definition.

\begin{definition}
Let $A$ be an affine bundle over $M$ of rank zero. A $k$-vector
$\bar{P}_{\hat{A}-\{O\}}$ on $\hat{A}-\{O\}=M\times (\R-\{0\})$ is
said to be homogeneous (sometimes called also linear) if
\[
\bar{P}_{\hat{A}-\{O\}}(d(tf_1),\dots,
d(tf_k))=(th)_{|\hat{A}-\{0\}},
\]
for $f_1,\dots ,f_k\in C^\infty(M,\R)$, where $t$ is the usual
coordinate on $\R$ and $h\in C^\infty(M,\R)$. \end{definition}

 Now, we deduce

\begin{corollary}\label{c3.6}
Let  $\pi:A\to M$ be an affine bundle on $M$ of rank $n$. Then:
\begin{enumerate}
\item If $n>0,$ there is a one-to-one  correspondence  between
affine Jacobi brackets $\{\cdot,\cdot\}_{(\zL_A,E_A)}$ on $A$ and
homogeneous Poisson brackets $\{\cdot,\cdot\}_{\bar\zL_{\hat A}}$
on the vector bundle $\hat A$, uniquely determined by the equation
\[
\{ a ,b \}_{(\zL_A,E_A)}=(\{ \bar{a},\bar{b}\}_{\bar{\zL}_{\hat
A}})_ {\vert A}
\] for $a,b:A\to \R$ affine functions  on $A$. The Jacobi structure
$(\zL_A,E_A)$ on $A$ is the restriction to $A$ of the Jacobi
structure $(\bar{\zL}_{\hat{A}}-\zD_{\hat{A}}\we
E_{\hat{A}},E_{\hat{A}})$ on $\hat{A},$ where $E_{\hat{A}}$ is the
hamiltonian vector field of the linear function $\zi_{\tilde
1}:\hat{A}\to \R$ with respect to $\bar\zL_{\hat{A}}$.\label{afpo}
\item
If $n=0,$ there is a one-to-one  correspondence  between  affine
Jacobi brackets $\{\cdot,\cdot\}_{(\zL_A,E_A)}$ on $A$ and linear
Poisson tensors ${\bar\zL_{\hat A-\{O\}}}$ on  $\hat
A-\{O\}=M\times (\R-\{0\})$, uniquely determined by the equation
\[
\{tf,tg\}_{\bar{\zL}_{\hat
A-\{O\}}}=\bar{\Lambda}_{\hat{A}-\{O\}}(d(tf),
d(tg))=t(\{f,g\}_{(\Lambda,E)})_{|\hat{A}-\{O\}}
\] for all $f,g\in C^\infty(M\times (\R-\{0\},\R)),$ where
$t$ is the usual coordinate on $\R$ and $(\Lambda,E)$ is the
Jacobi structure on $M$ induced by $(\zL_A,E_A).$ The linear
Poisson tensor ${\bar\zL_{\hat A-\{O\}}}$ is given by
\[
{\bar\zL_{\hat A-\{O\}}}=\frac{1}{t}\Lambda + \frac{\partial
}{\partial t}\wedge E.
\]
\end{enumerate}
\end{corollary}

\begin{remark}
{\rm Let $\pi:A\to M$ be an affine bundle on $M$ of rank $n$,
$n>0,$ and $\hat{A}$ the dual space of $A^+$ . If $({\cal
T}A=TA\oplus \R,(0,1))$ is the Jacobi algebroid of first-order
differential operators on $A$, then $P+ I\wedge Q\in
\Gamma(\wedge^k {\cal T}A)$ is affine if $(P + I\wedge Q)
(a_1,\dots, a_k)$ is an affine function, for all $a_1,\dots
,a_k:A\to \R$ affine functions on $A$.

On the other hand,  a $k$-vector $\bar{P}\in
\Gamma({\wedge^kT\hat{A}})$ on $\hat{A}$ is linear if
$\bar{P}(\bar{a}_1, \dots, \bar{a_k}$) is a linear function, for
all $\bar{a}_1, \dots ,\bar{a}_k$ linear functions on $\hat{A}$.

Now, let $\bar{P}\in \Gamma(\wedge^kT\hat{A})$ be a linear
$k$-vector on $\hat{A}.$ We consider the $k$-section $P'$ and the
$(k-1)$-section $Q'$ on $\hat{A}$ given by
\[
P'=\bar{P} - \Delta_{\hat{A}}\wedge
i(d\iota_{\tilde{1}})\bar{P},\makebox[1cm]{}
Q'=i(d\iota_{\tilde{1}})\bar{P}.
\]
Then, the restrictions $P$ and $Q$ to $A$ of $P'$ and $Q',$
respectively,  are tangent to $A$ and $P + I\wedge Q \in
\Gamma(\wedge^k{\cal T}A)$ defines an affine first-order
differential operator  on $A$. In fact, we have that this
correspondence between linear $k$-vectors on $\hat{A}$ and affine
first-order $k$-differential operators is one-to-one and that
\[ (P + I\wedge Q)(a_1,\dots
,a_k)=\bar{P}(\bar{a}_1,\dots ,\bar{a}_k))_{|A},
\]
for all $a_1,\dots ,a_k$ affine functions on $A.$ Here
$\bar{a}_i:\hat{A}\to \R$ denotes the linear function associated
with $a_i:A\to \R$. Moreover,  if $\bar{P_1}$ (respectively,
$\bar{P_2}$) is a linear $k_1$-vector (respectively, $k_2$-vector)
on $\hat{A}$ and $[\cdot,\cdot]_1$ is a Jacobi-Schouten bracket on
$A$  (see Section \ref{seccion2}) then
\[
[P_1+I\wedge Q_1,P_2 + I\wedge Q_2]_1=[\bar{P_1},\bar{P_2}]_{SN},
\] where $P_1+ I\wedge Q_1$ and $P_2+ I\wedge Q_2$ are the corresponding affine
first-order differential operators associated with $\bar{P_1}$ and
$\bar{P_2}$, respectively. The details and proofs of these results
can be found in \cite{last}.

Using the above facts one may directly deduce the first part of
Corollary \ref{c3.6}.

In the case $n=0$, if $\bar{P}_{\hat{A}-\{O\}}$ is a linear
$k$-vector on ${\hat{A}-\{O\}}=M\times (\R-\{0\})$ then we can
consider the $k$-vector $P'$ and the $(k-1)$-vector $Q'$ on
${\hat{A}-\{O\}}$ given by
\[
P'=t\bar{P}_{\hat{A}-\{O\}}-t\frac{\partial}{\partial t}\wedge
i(dt)\bar{P}_{\hat{A}-\{O\}}\;\;\;\;
Q'=i(dt)\bar{P}_{\hat{A}-\{O\}}.\] The restrictions $P$ and $Q$ to
$M$ of $P'$ and $Q'$, respectively, are tangent to $M$ and $P +
I\wedge Q\in \Gamma(\wedge^k(TM\oplus \R))$ defines an affine
first-order differential operator on $M.$ Moreover, we have that
this correspondence between linear $k$-vectors on $\hat{A}-\{0\}$
and affine first-order $k$-differential operators is bijective.
Note that the relation between local Lie algebras on rank 1 vector
bundles $L$ and homogeneous Poisson brackets on $L-\{ 0\}$ has
been first established by C.-M.~Marle \cite{Ma}.}\end{remark}

\section{Examples}\label{Examples} \setcounter{equation}{0} In
this section we present some examples and applications of the
above section.

${\bf 1.-}$   {\em Affine Poisson structures and special Lie
algebroid structures.} Let $(V,X)$ be a special vector bundle over
a manifold $M$. {\it A special Lie algebroid (resp. QD-algebroid)
structure} on $(V,X)$ is  a Lie algebroid (resp. QD-algebroid)
structure $(\lcf\cdot,\cdot\rcf, \rho)$ on $(V,X)$ for which the
section $X$ belongs to the center of the Lie algebra
$(\Gamma(V),\lcf\cdot,\cdot\rcf)$, that is, $\lcf X,Y\rcf=0,$ for
all $Y\in\Gamma(V).$

Then, if $A$ is an affine bundle with rank non-zero (resp. zero),
one can deduce from Theorems \ref{p1} and \ref{p2} and Corollary
\ref{c1} that there is a one-to-one correspondence between affine
Poisson structures on $A$ and special Lie algebroid (resp.
QD-algebroid)  structures on $A^\dag=(A^+,\tilde 1)$.

\medskip
${\bf 2.-}$ {\em Affine Jacobi structures on an affine space and
Lie algebra structures.} Let $A$ be an affine space of finite
dimension $n>0$  modeled on the space vector $V$. Then, using
Corollary \ref{c1}, we deduce that there is a one-to-one
correspondence between affine Jacobi structures on $A$ and Lie
algebra structures on $A^+$.

In the particular case, when $A$ is a vector space $V$, we have a
one-to-one correspondence between affine Jacobi structures on $V$
and Lie algebra structures on $V^*\times \R$, $V^*$ being the dual
vector space of $V$.

As a consequence of these facts and Example $1$, we obtain a
well-known result (see \cite{Bh}) which establishes a bijection
between affine Poisson structures on the vector space $V$ and
central extensions of Lie algebra structures on $V^*.$

\medskip
${\bf 3.-}$ {\em Affine Jacobi structures and triangular
generalized Lie bialgebroids.} We recall that a {\it triangular
generalized Lie bialgebroid} is a triple $((A,\lcf
\cdot,\cdot\rcf,\rho),\phi_0,P)$, where $A$ is a vector bundle
over $M$, $(\lcf\cdot,\cdot\rcf,\rho)$ is a Lie algebroid
structure on $A$, $\phi_0\in \Gamma(A^*)$ is a $1$-cocycle and
$P\in \Gamma(\wedge^2 A)$ a bisection on $A$ satisfying $\lcf
P,P\rcf+2P\wedge i(\phi_0)P=0$ (see \cite{IM2}).

Assume  that $\phi_0$ is nowhere vanishing and consider the affine
bundle $A_{\phi_0}=(A^*,\phi_0)^\ddag.$ A direct computation
proves that the (special) vector dual $A_{\phi_0}^+$ of
$A_{\phi_0}$ is isomorphic to the dual bundle $A^*$ of $A$. Thus,
the vector bundles $\hat{A}_{\phi_0}$ and $A$ are isomorphic.

Now, we consider the Poisson complete lift to $A\cong
\hat{A}_{\phi_0}$ of $P$ and $\phi_0$ given by
\[
\hat  P^c_{\phi_0}=P^c-\zi_{\phi_0} P^v+\zD_A\we(i_{\phi_0}P)^v.
\]
$\hat{P}_{\phi_0}^c$ is a linear Poisson structure on $A\cong
\hat{A}_{\phi_0}$ (see \cite{GM}). Next, we will see which is the
affine Jacobi structure on $A_{\phi_0}$ induced by
$\hat{P}_{\phi_0}^c.$

 In view of Corollary \ref{c3.6}, this
structure is the restriction to $A_{\phi_0}$   of the Jacobi
structure $(\hat P^c_{\phi_0}-\zD_A\we E_{\phi_0},E_{\phi_0})$,
where $E_{\phi_0}$ is  the  hamiltonian vector  field  of
$\zi_{\phi_0}$. Using  the identities $\zD_A(\zi_\zm)=\zi_\zm$,
$i_{\D\zi_\zm}X^c=(i_\zm X)^c+\zi_{i_X(d\mu)}$  and
$i_{\D\zi_\zm}X^v=(i_\zm X)^v$ which are valid for any 1-form
$\zm$, we get \[\begin{array}{rcl}
E_{\phi_0}&=&i_{\D\zi_{\phi_0}}(P^c-\zi_{\phi_0}
P^v+\zD_A\we(i_{\phi_0}P)^v)\\
&=&(i_{\phi_0}P)^c-\zi_{\phi_0}(i_{\phi_0}
P)^v+\zi_{\phi_0}(i_{\phi_0} P)^v-
(i_{{\phi_0}}i_{{\phi_0}}P)^v\zD_A=(i_{\phi_0} P)^c.
\end{array}\] Hence
\[
\hat  P^c_{\phi_0}-\zD_A\we  E_{\phi_0}= P^c-\zi_{\phi_0}
P^v-\zD_A\we((i_{\phi_0} P)^c-(i_{\phi_0}P)^v), \] and the Jacobi
structure on $A_{\phi_0}$ is the restriction of the Jacobi
structure
\[
(P^c-\zi_{\phi_0}P^v-\zD_A\we((i_{\phi_0}P)^c-(i_{\phi_0}P)^v),(i_{\phi_0}
P)^c). \] Note that $\zi_{\phi_0}=1$ on $A_{\phi_0}$. Let now
$I_0$ be  a  section  of  $A_{\phi_0}$.  We have the
decomposition $A=A_0\oplus\langle I_0\rangle\simeq A_0\oplus\R$,
where $A_0=V(A_{\phi_0})=Ker({\phi_0})$ is a 1-codimensional
vector subbundle of $A$. Since $d{\phi_0}=0$, $A_0$ is a  Lie
subalgebroid. Using the canonical linear coordinate $s$ in the
1-dimensional subbundle $\langle I_0\rangle\simeq\R$, we have that
\[
(I_0)^v=\pa_s,\;\;\;\;\; \zi_{\phi_0}=s,\;\;\;\;\; Q_0^c=Q_0^{c_0}
+ s\lcf I_0,Q_0\rcf^{v_0},
\]
for $Q_0\in \Gamma(\wedge^pA_0)$, where $c_0$ and $v_0$ denote the
complete and vertical lift of the Lie algebroid $A_0$. Here, of
course, we understand tensors on $A_0$ as tensors on $A\simeq
A_0\ti\R$ in obvious way. Note that if we identify $A_{\phi_0}$
with $A_0$ via the translation by $I_0$, then the restriction of
$Q_0^c$ to $A_{\phi_0}$ is tangent to $A_{\phi_0}$ and such a
restriction is the complete lift of $Q_0$ with respect to the {\it
Lie affgebroid } structure on $A_{\phi_0}$ in the terminology of
\cite{GGU}. A Lie affgebroid is a possible generalization of the
notion of a Lie algebroid to affine bundles. The main motivation
of the study of this concept was to create a geometrical model
which would be a natural environment for a time-dependent version
of Lagrange equations on Lie algebroids  (cf.
\cite{GGU,GGU1,MMS,SMM}).

Now, we can decompose $P=\zL + I_0\wedge E,$ where $\zL\in
\Gamma(\wedge^2A_0)$ and $E\in \Gamma(A_0),$ and we deduce
\[
P^c=\zL^{c_0} + s\lcf I_0,\zL\rcf^{v_0} + I_0^c \wedge E^{v_0} +
\partial_s\wedge (E^{c_0} + s\lcf I_0, E\rcf^{v_0})\makebox[2cm]{and}
P^v=\zL^{v_0} + \partial_s\wedge E^{v_0}.
\]
Thus, writing $\zD_A=\zD_{A_0}+s\pa_s$, we get finally the Jacobi
structure on $A_{\phi_0}$ (identified with $A_0$ via the
translation by $I_0$) in the form
\[
(\zL^{c_0}-\zL^{v_0} + \lcf I_0,\zL\rcf^{v_0}+(I_0)^c\we
E^{v_0}-\zD_{A_0}\we(E^{c_0}-E^{v_0} + \lcf
I_0,E\rcf^{v_0}),E^{c_0}). \] In particular, if $\rho(I_0)=0$ and
$I_0$ is central, i.e., $(A,I_0)$ is a special Lie algebroid, then
$(I_0)^c=0$, so we end up with the Jacobi structure \be\label{eee}
(\zL^{c_0}-\zL^{v_0} -\zD_{A_0}\we(E^{c_0}-E^{v_0}),E^{c_0}). \ee

Now, we consider a particular example of triangular generalized
Lie bialgebroid.

 Let $(A,\lcf\cdot,\cdot\rcf_A,\rho_A)$ be a Lie
algebroid over $M$ and let $(\Lambda,E)\in \Gamma(\wedge^2A)\oplus
\Gamma(A)$ be a pair satisfying the following properties
\[
\lcf\Lambda,\Lambda\rcf_A=-2\Lambda\wedge E,\;\;\;\;\;\;
\lcf\Lambda,E\rcf_A=0.
\]
Here $\lcf\cdot,\cdot\rcf_A$ denotes the Schouten bracket
associated with the Lie algebroid $A.$

We will prove that, in such a case, it is possible to define an
affine Jacobi structure over $A$. In fact, we can consider the Lie
algebroid     structure $(\lcf\cdot,\cdot\rcf_{A_1},\rho_{A_1})$
over $A_1=A\oplus\R$ given as follows
\[
\begin{array}{rcl}
\lcf(X,f),(Y,g)\rcf_{A_1}=(\lcf
X,Y\rcf_A,\rho_A(X)(g)-\rho_A(Y)(f)),\;\;\;\;
\rho_{A_1}(X,f)=\rho_A(X),
\end{array}
\]
for all $(X,f),(Y,g)\in \Gamma(A)\times
C^\infty(M,\R)\cong\Gamma(A_1).$ The pair $(0,1)\in
\Gamma(A^*)\times C^\infty(M,\R)\cong \Gamma(A_1^*)\cong
\Gamma(A^*\oplus\R)$ defines a $1$-cocycle for this algebroid.
Moreover, since $P=\zL+I\we E$ is a  canonical  structure  for the
corresponding Schouten-Jacobi bracket, $(A_1, (0,1), (\Lambda,E))$
is a triangular generalized Lie bialgebroid (see \cite{IM2}). In
this case, $I_0=(0,1)\in \Gamma(A_1)$ is central and
$\rho_{A_1}(I_0)=0.$ Thus, we have defined an affine Jacobi
structure on $A$ given by \be\label{aaa}
(\zL_A,E_A)=(\zL^{c}-\zL^v -\zD_{A}\we(E^{c}-E^v),E^c), \ee where
the complete lifts are lifts for the Lie algebroid $A$.

 Now, if we consider the
particular case when
$(A,\lcf\cdot,\cdot\rcf_A,\rho_A)=(TM,[\cdot,\cdot]_{SN},1_{TM})$
and $(M,\Lambda,E)$ is a Jacobi manifold, then we have that  the
affine Jacobi structure $(\Lambda_{TM},E_{TM})$ on $TM$ ({\it the
affine tangent Jacobi structure on $TM$}) is given by (see
(\ref{aaa}))
\[
\Lambda_{TM}=\Lambda^c-\Lambda^v-\Delta_{TM}\wedge (E^c-E^v),
\;\;\;\; E_{TM}=E^c,
\]
where $\Delta_{TM}$ is the Liouville vector field on $TM,$
$\Lambda^c$ (resp. $E^c$) and $\Lambda^v$ (resp. $E^v$) are the
complete and vertical  lift of $\Lambda$ (resp. $E$). This
structure was first considered by Vaisman in \cite{V2}. If $E=0$
(that is, $(M,\zL)$ is a Poisson manifold) we obtain {\it the
affine tangent Poisson structure} on $TM$ given by $
\Lambda_{TM}=\Lambda^c-\Lambda^v.$

Note that the linear  tangent Poisson structure on $TM$ is
$\Lambda^c$.

${\bf 4.-}$ {\em Homogeneous Jacobi structures and Jacobi
algebroids.} A Jacobi structure $(\zL,E)$ on a manifold $M$  is
called {\it homogeneous of degree $k$} with respect to  a vector
field $\zD$ on $M$  if $\zL$ and $E$ are homogeneous of degree $k$
with respect to $\zD$. A homogeneous Jacobi structure of degree
$-1$ we will just call {\it homogeneous}. By a  {\it homogeneous}
Jacobi structure on a vector bundle  $A$,  we will always
understand a Jacobi structure which is homogeneous with respect to
the Liouville vector field $\zD_A$. Then, we have the following
characterizations.

\begin{theorem}\label{ho} Let $(\zL,E)$ be a Jacobi structure on a  vector
bundle $A$. Then, the following are equivalent:
\begin{description}
\item{(a)} $(\zL,E)$ is homogeneous;
\item{(b)} The Jacobi bracket $\{\cdot,\cdot\}_{(\zL,E)}$ is  linear
and affine and  the bracket of a linear function and the constant
function 1 is  a  basic function;
\end{description}
\end{theorem}
{\bf Proof.-} $(a)\Rightarrow (b)$ If $\zm,\zn$ are sections of
$A^*$ then, from (\ref{lineal}), it follows that
\[
D(\{\zi_\zm,\zi_\zn\}_{(\Lambda,E)})=
\{D(\zi_\zm),\zi_\zn\}_{(\Lambda,E)}+\{\zi_\zm,D(\zi_\zn)\}_{(\Lambda,E)}=
0,
\]
for $D=\zD_A-I$, which implies that
$\{\zi_\zm,\zi_\zn\}_{(\Lambda,E)}$ is linear.

On the other hand, since $E$ is homogeneous, we obtain that $
\zD_A(E(\zi_\zm))=0$  and thus
$E(\zi_\zm)=\{1,\zi_\mu\}_{(\Lambda,E)}$ is a basic function.
Consequently, using the results of \cite{IM1}, we deduce that
$\{\cdot,\cdot\}_{(\Lambda,E)}$ is affine.

$(b)\Rightarrow (a)$ Let $f$ be a basic function. If $\zm$ is a
section of $A^*$ then $
\{1,f\zi_\zm\}_{(\Lambda,E)}=E(f\zi_\zm)=E(f)\zi_\zm +
f\{1,\zi_\zm\}_{(\Lambda,E)}$ is a basic function. Therefore,
$E(f)=0.$ Moreover,
$\Delta_A(E(\iota_\mu))=\Delta_A(\{1,\iota_\mu\}_{(\Lambda,E)})=0.$
Consequently, $[\Delta_A,E]=-E.$

Next, we will prove that $\zL$ is linear. For $\mu,\nu\in
\Gamma(A^*)$, we have
\[
\Lambda(\D\zi_\zm,\D\zi_\zn)=\{\zi_\zm,\zi_\zn\}_{(\Lambda,E)}-\zi_\zm
E(\zi_\zn) + \zi_\nu E(\zi_\zm)
\]
since $E(\zi_\zm)$ and $E(\zi_\zn)$ are basic functions, we
conclude that $\Lambda(d\iota_\mu,d\iota_\zn)$ is a linear
function. This implies that  $\Lambda$ is homogeneous (see Section
\ref{seccion2}). \hfill$\Box$

Let $\pi:A\to M$ be a vector bundle over the  manifold $M$ of rank
$n$, $n>0$. Consider a homogeneous Jacobi structure $(\Lambda_A,
E_A)$ on $A$. Then $(\Lambda_A,E_A)$ is an affine Jacobi structure
on $A$ (see Theorem \ref{ho}).

Thus, from Theorem \ref{p1}, we have that there exists a Lie
algebroid structure on $A^+=A^\ast \oplus\R$. It is not difficult
to show that such a Lie algebroid structure is given by
\[
\begin{array}{rcl}
\lcf(\alpha,f),(\beta,g)\rcf^{+}=(\lcf
\alpha,\beta\rcf_*,\rho_*(\alpha)(g)-\rho_*(\beta)(f)+(g\alpha-f\beta)(X_0)),\;\;\;\;\;
\rho^+(\alpha,f)=\rho_*(\alpha)
\end{array}
\]
for all $(\alpha,f), (\beta,g)\in \Gamma(A^*)\times
C^\infty(M,\R)\cong \Gamma(A^+)$, where $(\lcf\cdot,\cdot\rcf_*,
\rho_*)$ is a Lie algebroid structure on $A^*$ and $X_0$ is a
$1$-cocycle on $A^*$.

Using the above fact and Corollary \ref{c1}, we conclude that
there exists a one-to-one correspondence between homogeneous
Jacobi structures  on $A$ and Jacobi algebroid structures on
$A^*$. This result was proved in \cite{IM1}.

\section{Affine-homogeneous Jacobi structures on an affine
bundle}\label{strongly-affine} \setcounter{equation}{0}

Let $\pi:A\to M$ be a vector bundle over a manifold $M$. We stress
the existence of the canonical family ${\cal V}(A)=\{ X^v:X\in
\Gamma(A)\}$ of vertical lifts of sections of $A$. These are
pair-wise commuting vector fields. They can be viewed as invariant
vector fields on $A$ when we view $A$ as  a commutative Lie
groupoid. Any tensor $Y$  on  $A$  is called {\it invariant} if it
is invariant with respect to ${\cal V}(A)$, i.e., the Lie
derivative $\Ll_XY$ vanishes for every $X\in{\cal V}(A)$. It is
easy to see that the set  of invariant vector fields coincides
with  ${\cal V}(A)$.

It is worth noticing that linear tensor fields have a special
property closely related to the fact that they live on vector
bundles. On the zero-section of  a vector bundle $A$ we have the
full decomposition of $T_xA$ into the vertical  and the horizontal
parts. This, of course, makes sense for any contravariant tensor
and we will say that a $r$-vector $\zL$ on $A$ is {\it vertically
vanishing on the zero-section} if its vertical part vanishes on
the zero-section. This simply means that $
\Lambda(\D\zi_{\zm_1},\dots ,\D\zi_{\zm_r})\circ O=O,$  for
$\zm_1,\dots, \zm_r\in \Gamma(A^*)$, where $O$ is the zero-section
of $A$. In fact, if $\zL$ is a $r$-vector on $A$ then one may
define a section $V_\zL$ of the vector bundle $\wedge^rA\to M$ as
follows. If $\zm_1,\dots ,\zm_r\in \Gamma(A^*)$ then
\[
V_{\zL}(\zm_1,\dots ,\mu_r)=\zL(\D\zi_{\zm_1},\dots
,\D\zi_{\zm_r})\circ O.
\]
Now, a $2$-vector $\zL$ on $A$ is {\it affine} if $\zL(da,db)$ is
an affine function, for $a,b:A\to \R$ affine functions on $A$. It
is easy to prove that  $\zL$ is affine if and only if
$\zL-V_\zL^v$ is a linear $2$-vector on $A$, where $V_\zL^v$ is
the standard vertical lift of $V_\zL\in \Gamma(\wedge^2A).$
Furthermore, we have

\begin{lemma}\label{l2.0}
Let $\zL$ be a $2$-vector on a vector bundle $A$. Then, $\zL$ is
affine if and only if $\Ll_X\zL$ is an invariant bivector field
for any invariant vector field $X\in {\cal V}(A)$. Moreover, $\zL$
is linear if and only if it is affine and vertically vanishes on
the zero-section $A$.
\end{lemma}
{\bf Proof.-} The proof is obvious and depends on the fact that if
$f:A\to \R$ is a smooth real function on $A$ then $f$ is affine if
and only if $X(f)$ is a basic function, for any invariant vector
field $X\in {\cal V}(A)$. In addition, $f$ is linear if and only
if $X(f)$ is a basic function, for all $X\in {\cal V}(A)$ and
$f\circ O$ identically vanishes. Using these facts and
(\ref{basica}), we deduce the result. \hfill$\Box$

We recall that a $2$-vector $\Lambda$ on a vector bundle is linear
if and only if it is homogeneous (with respect to the Liouville
vector field of $A$). On the other hand, a Jacobi structure
$(\Lambda,E)$ on $A$  is homogeneous if $\Lambda$ and $E$ are
homogeneous (see Example $4$ in Section \ref{Examples}). Next,
using the definition of invariant tensor fields on $A$, we will
characterize homogeneous Jacobi structures.

\begin{theorem}
Let $(\Lambda,E)$ be a Jacobi structure on a vector bundle $A.$
Then, the following sentences  are equivalent:
\begin{enumerate}
\item[(a)] $(\Lambda,E)$ is homogeneous;
\item[(b)] $E\in{\cal V}(A)$  and   there is a linear  Poisson
structure   $\bar   \zL$    such    that
\[
\zL=\bar\zL+E\we\zD_A;
\]

\vspace{-15pt}

\item[(c)] $E\in{\cal V}(A)$, $\Ll_X\zL$ is an  invariant bivector field for
any invariant vector field $X\in {\cal V}(A)$ and $\zL$ vanishes
vertically on the zero-section of $A$.
\end{enumerate}
\end{theorem}
{\bf Proof.-}  $(a)\Rightarrow (b)$ Let $f$ be a basic function.
If $\mu$ is a section of $A^*$ then, from Theorem \ref{ho}, we
deduce that the function
$\{1,f\iota_\mu\}_{(\Lambda,E)}=E(f\iota_\mu)=E(f)\iota_\mu +
f\{1,\iota_\mu\}_{(\Lambda,E)}$ is a basic function. Thus,
$E(f)=0$. Therefore, since $E(i_{\zn})$ is a basic function, for
all $v\in \Gamma(A^*),$ we conclude that $E\in {\cal V}(A).$

Now, using that the Jacobi bracket $\{\cdot,\cdot\}_{(\Lambda,E)}$
is linear (see Theorem \ref{ho}), we obtain that $\Lambda$ is
linear. This implies that ${\Ll}_{\Delta_A}\Lambda=-\Lambda$ and,
since that $E\in {\cal V}(A)$, it follows that
$[\bar{\Lambda},\bar{\Lambda}]_{SN}=0$, that is, $\bar\Lambda$ is
a Poisson structure on $A.$ Finally, from (\ref{lineal}) and using
that $\Lambda$ is linear and the fact that $E\in {\cal V}(A)$, we
have that $\bar{\Lambda}$ is a linear $2$-vector on $A.$

$(b)\Rightarrow (c)$  It follows from Lemma \ref{l2.0}.

$(c)\Rightarrow (a)$ If $E\in {\cal V}(A),$ it is clear that $E$
is homogeneous. Therefore, using again Lemma \ref{l2.0}, we have
that $\Lambda$ is linear and, consequently, $\zL$ is homogeneous.

\hfill$\Box$

 Note that the description of homogeneity of tensor
fields given in  $(c)$ of  the  above  theorem  is  a  variant  of
the description   of   {\it multiplicative}  tensors  on  Lie
groups or  Lie  groupoids,  here  in the
commutative case.

We can also try to define homogeneity for affine  bundles. In
order to do this we  describe homogeneous Jacobi structures, using
the Schouten-Jacobi bracket $[\cdot,\cdot]_1$ on the Grassmann
algebra of  first-order polydifferential operators.

\begin{theorem}\label{ho1}
$(a)$ A Jacobi structure $(\zL,E)$ on a manifold $M$ is
homogeneous with respect to a vector field $\zD$ if and only if
\be
[\zD-I,\zL+I\we E]_1=0, \ee where $[\cdot,\cdot]_1$ is the
canonical Schouten-Jacobi bracket on the  Grassmann  algebra  of
first-order  polydifferential  operators on $M$.

\noindent $(b)$ A Jacobi structure  $(\zL,E)$  on  a  vector
bundle $A$ is homogeneous if and only if the bivector field $\zL$
vertically vanishes on the zero-section of $A$ and \be\label{r5.2}
[(X_1)^v,[(X_2)^v,\zL+I\we  E]_1]_1=0, \mbox{ for all
$X_1,X_2\in\zG(A)\oplus C^\infty(M,\R)$}.\ee
\end{theorem}
{\bf Proof.-} (a) It is a direct consequence from  $
[\zD-I,\zL+I\we E]_1=[\zD,\zL]_{SN}+I\we[\zD,E]_{SN}+\zL+I\we E$.

(b) Suppose that $\Lambda$ vertically vanishes on the zero-section
of $A$ and that (\ref{r5.2}) holds. Then, since for any
$X\in\zG(A)$, $ [X^v,[(1_M)^v,\zL+I\we E]_1]_1=[E,X^v]_{SN}=0, $
the vector field $E$ is invariant, i.e., $E\in{\cal V}(A)$. Hence,
for any $X_1,X_2\in\zG(A)$, $ [X_1^v,[X_2^v,\zL+I\we
E]_1]_1=[X_1^v,[X_2^v,\zL]_{SN}]_{SN}=0, $ i.e., $\zL$ is affine.
Thus, from Lemma  \ref{l2.0}, we deduce that $\Lambda$ is linear.

Conversely, if $\zL$ is linear and $E$ is invariant, then for any
$X\in\zG(A)$, $f\in C^\infty(M,\R)$, $$ [X^v+f^v,\zL+I\we
E]_1=[X^v,\zL]_{SN}+[f^v,\zL]_{SN}-f^vE. $$ But $[X^v,\zL]_{SN}$
and $([f^v,\zL]_{SN}-f^vE)$ are invariant tensors and the theorem
follows. \hfill$\Box$

Since on an affine bundle $A$ no Liouville vector field and no
linear functions exist, we will use a concept of homogeneity
suggested by Theorem \ref{ho1}. The concept  of invariance  is
clear:  the model vector bundle $V(A)$ acts  on $A$  by
translations,  so  we can  lift vertically sections $X$ of $V(A)$
to vector fields $X^v$ on $A$. We lift vertically functions on $M$
to functions on $A$ in obvious  way.  Denote the vector  space  of
first-order differential  operators  spanned by vertical lifts of
both types by ${\cal V}_1(A)$. It is easy to see that ${\cal
V}_1(A)$ is a maximal subalgebra in the Lie algebra of all  first
order differential operators on $A$. A first-order
polydifferential operator $F$  on $A$ is called {\it
affine-invariant} if it is invariant with respect to elements of
${\cal V}_1(A)$, i.e., $[D,F]_1=0$ for any $D\in{\cal V}_1(A)$,
where $[\cdot,\cdot]_1$ is the canonical Schouten-Jacobi bracket
on the Grassmann algebra of first-order  polydifferential
operators on $A$ and {\it affine-homogeneous} if $[D,F]_1$ is
affine-invariant for any $D\in{\cal V}_1(A)$. In other words, $F$
is affine-homogeneous if $[D_1,[D_2,F]_1]_1=0$     for all
$D_1,D_2\in{\cal V}_1(A)$. In particular, we propose the following

\begin{definition}\label{d5.1''} A Jacobi structure $(\zL_A,E_A)$ on the
affine bundle $A\to M$ is said to be affine-ho\-mo\-ge\-ne\-ous if
$[D_1,[D_2,\zL_A+I\we E_A]_1]_1=0$     for all $D_1,D_2\in{\cal
V}_1(A)$.
\end{definition}

Although the Jacobi bracket associated with an affine Jacobi
manifold $(A,\Lambda,E)$ is closed with respect to the affine
functions, the hamiltonian vector field $X_a^{(\Lambda,E)}$ of an
affine function $a:A\to \R$ is not, in general, affine, that is,
if $b:A\to \R$ is an affine function then $X^{(\Lambda,E)}_a(b)$
is not, in general,  an affine function.

In fact, if $A$ is the vector space $\R^3$ and we consider the
Jacobi structure $(\Lambda,E)$ on $A$ given by

\begin{equation}\label{ejemplo}
\Lambda=x_1x_3\frac{\partial }{\partial x_2}\wedge \frac{\partial
}{\partial x_3} - x_1^2\frac{\partial }{\partial x_1}\wedge
\frac{\partial }{\partial x_2}, \;\;\;\; E=x_1\frac{\partial
}{\partial x_2},
\end{equation}
where $(x_1,x_2,x_3)$ are canonical coordinates in $\R^3$, then
$(\Lambda,E)$ is an affine Jacobi structure over $\R^3$ and
$X_{x_2}^{(\Lambda,E)}(x_3)=x_1x_3$ is not an affine function.
Note that, in the case of affine Poisson structures, the
hamiltonian vector fields of affine functions are affine.

\begin{definition} An affine Jacobi structure $(\Lambda_A,E_A)$  over an
affine bundle $\pi:A\to M$ is a strongly-affine Jacobi structure
if  the hamiltonian vector field $X^{(\Lambda_A,E_A)}_a$ is
affine, for every affine function $a:A\to \R$.
\end{definition}

The following result relates affine-homogeneous and
strongly-affine Jacobi structures.

\begin{proposition}\label{strong}
Let $(\Lambda_A,E_A)$ be a  Jacobi structure over an affine bundle
$\pi:A\to M$ of rank $n$, $n>0$. Then, the following sentences are
equivalent:
\begin{description}
\item{(i)} $(\Lambda_A,E_A)$ is affine-homogeneous;
\item{(ii)} $E_A$ is affine-invariant and $\zL_A$ is affine;
\item{(iii)} $(\Lambda_A,E_A)$ is strongly-affine;
\item{(iv)} $(\Lambda_A,E_A)$ is affine and  basic  functions  form  an
ideal  in  the  algebra  of  affine  functions  with   respect to
the corresponding Jacobi bracket;
\item{(v)} $(\Lambda_A,E_A)$ is affine and there exists
$\bar{X}_0\in \Gamma(\hat{A})$ such that
\begin{equation}\label{EAf}
E_A(a)=-\bar{X}_0(\tilde{a})\circ \pi
\end{equation}
for all affine functions $a:A\to \R$, where $\tilde{a}\in
\Gamma(A^+)$ is the section of $A^+$ associated with the affine
function $a$.
\end{description}
\end{proposition}
{\bf Proof.-} $(i)\Leftrightarrow (ii)$ Proceeding as in the proof
of Theorem \ref{ho1} $(b)$ we deduce that the sentences $(i)$ and
$(ii)$ are equivalent.

$(ii)\Rightarrow (iii)$ Let $a$ be an affine function and consider
the hamiltonian vector field $X_a^{(\zL_A,E_A)}$ of $a$ with
respect to $(\zL_A,E_A)$. If $b:A\to \R$ is an affine function
then
\[
X_a^{(\zL_A,E_A)}(b)=\Lambda_A(da,db) + a E_A(b).
\]
Now, the condition $[X,E_A]=0,$ for all $X\in {\cal V}(A)$,
implies that $E_A\in {\cal V}(A)$ and, thus, $E_A(b)$ is a basic
function. Therefore, $X_a^{(\zL_A,E_A)}(b)$ is an affine function.

$(iii)\Rightarrow (iv)$ Let $a:A\to \R$ be an affine function. We
will prove that $E_A(a)=\{1,a\}_{(\zL_A,E_A)}$ is a basic
function. In fact, if $p$ is a point of $M$, we will show that
$(E_A(a)_{|A_p})$ is constant. For this purpose, we consider an
affine function $b:A\to \R$ such that the restriction to $V(A)_p,$
$\hat{b}_{|V(A)_p},$ of the linear map $\hat{b}$ associated with
$b$ is not zero. Then, we have that $
X^{(\Lambda_A,E_A)}_a(b)=\{a,b\}_{(\Lambda_A,E_A)} - bE_A(a)$ and
$\{a,b\}_{(\Lambda_A,E_A)} $ are  affine functions. Thus,
$bE_A(a)$ is an affine function and $E_A(a)_{|A_p}$ is constant.

Now, suppose that $f$ is a basic function. Then,
\[\{a,fb\}_{(\zL_A,E_A)}=f\{a,b\}_{(\zL_A,E_A)} +
b\{a,f\}_{(\zL_A,E_A)}-fb\{a,1\}_{(\zL_A,E_A)}
\]
is an affine function. This implies that $b\{a,f\}_{(\zL_A,E_A)}$
is an affine function and, therefore, $(\kern-1pt\{a,\kern-1pt
f\}_{(\zL_A,E_A)}\kern-1pt )_{|A_p}$ is constant. Consequently, we
have proved that $\{a,f\}_{(\zL_A,E_A)}$ is a basic function.

$(iv)\Rightarrow (v)$ Let $a:A\to \R$ be an affine function on
$A$. Then, the function $E_A(a)=\{1,a\}_{(\zL_A,E_A)}$ is basic.
Thus, there exists $\bar{X}_0(\tilde{a})\in C^\infty(M,\R)$ such
that
\begin{equation}\label{e5.2'}
E_A(a)=-\bar{X_0}(\tilde{a})\circ \pi.
\end{equation}
Next, we will prove that if $f$ is a basic function then
\begin{equation}\label{e5.2''}
E_A(f)=\{1,f\}_{(\zL_A,E_A)}=0.
\end{equation}
Suppose that $p$ is a point of $M$ and that $a:A\to \R$ is an
affine function such that the restriction to $V(A)_p$,
$\hat{a}_{|V(A)_p},$ of the linear map $\hat{a}$ associated with
$a$ is not zero. Then, $
\{1,fa\}_{(\zL_A,E_A)}=\{1,f\}_{(\zL_A,E_A)}a +
f\{1,a\}_{(\zL_A,E_A)}$ is a basic function and therefore,
$(\{1,f\}_{(\zL_A,E_A)})_{|A_p}=0.$ Consequently, from
(\ref{e5.2'}) and (\ref{e5.2''}), we deduce that
$\bar{X_0}:\Gamma(A^+)\to C^\infty(M,\R)$ is a section of the
vector bundle $\hat{A}=(A^+)^*.$

$(v)\Rightarrow (ii)$ Let $f$ be a basic function on $A$. Using
(\ref{EAf}), it follows that $aE_A(f)=0,$ for any affine function
$a:A\to \R$ on $A$. This implies that
\begin{equation}\label{e5.2'''}
E_A(f)=0.
\end{equation}
Consequently, from (\ref{EAf}) and (\ref{e5.2'''}), we obtain that
$E_A$ is affine-invariant.

On the other hand, if $a,b:A\to \R$ are affine functions then,$
\{a,b\}_{(\zL_A,E_A)}=\zL_A(da,db)+ aE_A(b)-bE_A(a)$ is an affine
function and, thus, $\zL_A(da,db)$ is again an affine function.
This proves that $\zL_A$ is affine. \hfill$\Box$
\begin{remark}
{\rm If $(\Lambda_A,E_A)$ is an affine-homogeneous Jacobi
structure on an affine bundle $A$ then, from Proposition
\ref{strong}, we can deduce that it is an affine Jacobi structure
and that the local expressions of $\Lambda_A$ and $E_A$ are as in
(\ref{e5}). Moreover, in this case, $c_{0\alpha}^\beta=0$, for all
$\alpha,\beta=1,\dots ,n$ .}
\end{remark}

Next, we will establish  a one-to-one correspondence between
strongly-affine Jacobi structures over an affine bundle $\pi:A\to
M$ and a particular class  of Lie algebroid structures on $A^+$.

\begin{definition} An {\it almost-special Lie algebroid structure} on
a special vector bundle $(V,X)$ is a Lie algebroid structure
$(\lcf\cdot,\cdot\rcf,\rho)$ on $V$ such that the submodule of
$\Gamma(V)$ generated by the section $X$ is an ideal of the Lie
algebra $(\Gamma(V),\lcf\cdot,\cdot\rcf).$
\end{definition}

\begin{remark}{\rm
 \begin{description}
\item[ $i)$] Every special vector bundle is of the  form
$(A^+,\tilde 1)$ for an affine bundle $A$.
\item[$ii)$]
A Lie algebroid structure $(\lcf\cdot,\cdot\rcf^{+},\rho^+)$ on
$A^+$ is almost-special if there exists a map
$\bar{X}_0:\Gamma(A^+)\to C^\infty(M,\R)$ such that
\begin{equation}\label{special} \lcf
\tilde{1},\tilde{a}\rcf^{+}=-\bar{X}_0(\tilde{a})\tilde{1},\;\;\;\mbox{for
all $\tilde{a}\in \Gamma(A^+).$}
\end{equation}
\end{description}}\end{remark}

 This type of Lie algebroids has the following properties.
\begin{proposition}\label{almost}
Let $\pi:A\to M$ be an affine bundle of rank $>0$, and
$(\lcf\cdot,\cdot\rcf^{+},\rho^{+} )$ be an almost-special Lie
algebroid structure  over $A^+.$ Then:
\begin{enumerate}
\item
$\rho^{+}(\tilde{1})=0,$
\item
$\bar{X}_0$ is a $1$-cocycle of the Lie algebroid
$(A^+,\lcf\cdot,\cdot\rcf^{+},\rho^+)$.
\end{enumerate}
\end{proposition}
{\bf Proof.-} Consider ${f_M}\in C^\infty(M,\R)$  and $p\in M.$ We
will prove that $\rho^{+}(\tilde{1})({f_M})(p)=0.$ In fact, if
$\tilde{a}\in \Gamma(A^+)$, then
\[
\bar{X_0}({f_M}\tilde{a})\tilde{1}=-\lcf
\tilde{1},{f_M}\tilde{a}\rcf^{+}=-{f_M}\lcf
\tilde{1},\tilde{a}\rcf^{+}
-
\rho^{+}(\tilde{1})({f_M})\tilde{a}={f_M}\bar{X}_0(\tilde{a})\tilde{1}
- \rho^{+}(\tilde{1})({f_M})\tilde{a}. \] Thus,
\begin{equation}\label{almost1}
({f_M}\bar{X}_0(\tilde{a})-\bar{X}_0
({f_M}\tilde{a}))\tilde{1}-\rho^{+}(\tilde{1})({f_M})\tilde{a}=0.
\end{equation}
If we consider $a:A\to \R$ an affine function such that the
associated linear function $\hat{a}$ satisfies $\hat{a}(p)\not=0$,
from (\ref{almost1}), one can deduce that
$\rho^{+}(\tilde{1})({f_M})(p)=0.$ So, we have $(i)$.

Substituting $(i)$ in (\ref{almost1}), we prove that $\bar{X_0}$
is $C^\infty(M,\R)$-linear.

Finally, using the Jacobi identity of $\lcf\cdot,\cdot\rcf^{+}$
and (\ref{almost1}), we obtain
\[
\begin{array}{lcl}
0=\lcf \tilde{1},\lcf\tilde{a},\tilde{b}\rcf^{+}\rcf^{+} - \lcf
\tilde{a},\lcf\tilde{b},\tilde{1}\rcf^{+}\rcf^{+}-\lcf
\tilde{b},\lcf\tilde{1},\tilde{a}\rcf^{+}\rcf^{+}
&=&-\bar{X}_0(\lcf\tilde{a},\tilde{b}\rcf^{+})\tilde{1}-
\lcf\tilde{a},\bar{X}_0(\tilde{b})\tilde{1}\rcf^{+} +
\lcf\tilde{b},\bar{X}_0(\tilde{a})\tilde{1}\rcf^{+}\\[5pt]
&=&(-\bar{X}_0(\lcf\tilde{a},\tilde{b}\rcf^{+})
+\rho^{+}(\tilde{a})(\bar{X}_0(\tilde{b}))-
\rho^{+}(\tilde{b})(\bar{X}_0(\tilde{a})))\tilde{1}.
\end{array}
\]
Therefore, $\bar{X}_0$ is a $1$-cocycle in $A^+.$ \hfill$\Box$

Now, as a consequence  of Corollary \ref{c1} and Propositions
\ref{strong} and \ref{almost}, we conclude that

\begin{corollary}\label{strong-almost}
There exists a one-to-one correspondence between strongly-affine
Jacobi structures over an affine bundle $\pi:A\to M$ of rank $n>0$
and almost-special Lie algebroid structures on $A^+$.
\end{corollary}

Now, suppose that $\pi:A\to M$ is a vector bundle and that
$A_1=A\times \R$. Let $(\Lambda_A,E_A)$ be a strongly-affine
Jacobi structure over $A$ and
$((\lcf\cdot,\cdot\rcf^{+},\rho^{+}),\bar{X}_0)$ be the associated
almost-special Lie algebroid structure on $A^+$. Under the
identification between $A^+$ and  $A_1^*$, the section $\tilde{1}$
of $A^+$ is the pair $(0,1)\in \Gamma(A^*)\times
C^\infty(M,\R)\cong \Gamma(A_1^*).$ Moreover, since
$\bar{X}_0(\tilde{1})=0$ and $\rho^{+}(\tilde{1})=0$ (see
(\ref{special}) and Proposition \ref{almost}), we deduce that
there exist maps
\[
\begin{array}{rcl}
\lcf\cdot,\cdot\rcf_*:\Gamma(A^*)\times \Gamma(A^*)\to
\Gamma(A^*),&&\rho_*:\Gamma(A^*)\to \X(M),\\ X_{0}:\Gamma(A^*)\to
C^\infty(M,\R),&& P_0:\Gamma(A^*)\times \Gamma(A^*)\to
C^\infty(M,\R)\end{array} \] such that
\begin{equation}\label{*}
\begin{array}{rcl}
\rho^{+}({a'},f_M)&=&\rho_*({a'})\\
\lcf({a'},{f}_M),({b'},{g_M})\rcf^{+}&=&(\lcf {a'},{b'}\rcf_*,
-P_0({a'},{b'})-{f_M}X_{0}({b'}) +
{g_M}X_{0}({a'})-\rho_*({b'})({f_M})\\&& + \rho_*({a'})({g_M})),
\end{array}
\end{equation}
for all $({a'},{f}_M), ({b'},{g}_M)\in \Gamma(A^*)\times
C^\infty(M,\R).$ A direct computation, using that
$(\lcf\cdot,\cdot\rcf^{+},\rho^{+})$ is a Lie algebroid structure,
shows that $(A^*,\lcf\cdot,\cdot\rcf_*,\rho_*)$ is a Lie
algebroid, that $X_{0}$ defines a $1$-cocycle in
$(A^*,\lcf\cdot,\cdot\rcf_*,\rho_*)$ and that
$P_0:\Gamma(A^*)\times \Gamma(A^*)\to C^\infty(M,{\R})$ is a
skew-symmetric $C^\infty(M,\R)$-bilinear mapping such that
$d_*P_0=-X_{0}\wedge P_0,$ where $d_*$ denotes the differential of
the  Lie algebroid $(A^*,\lcf\cdot,\cdot\rcf_*,\rho_*).$

\begin{remark}
{\rm Using that $X_0$ is a $1$-cocycle of the Lie algebroid
$(A^*,\lcf\cdot,\cdot\rcf_*,\rho_*)$ and that $d_*P_0=-X_0\wedge
P_0$, we deduce the following facts: \begin{description} \item[i)]
The map $\nabla:\Gamma(A^*)\times C^\infty(M,\R)\to
C^\infty(M,\R)$ given by
\[\nabla_{a'}f_M=\rho_*(a')(f_M) + a'(X_0)f_M,
\]
for $a'\in \Gamma(A^*)$ and $f_M\in C^\infty(M,\R)$ defines a
representation of the Lie algebroid $(A^\ast, \lcf \cdot ,\cdot
\rcf _*,\rho _*)$ on the vector bundle $M\times \R\to M$.

\item[ii)] $P_0:\Gamma(A^*)\times \Gamma(A^*)\to C^\infty(M,\R)$ is a
$2$-cocycle for the representation $\nabla$.
\end{description}

Thus, in the terminology of Mackenzie \cite{Mk} (see pag. 205-206
in \cite{Mk}), the Lie algebroid
$(A_1^*,\lcf\cdot,\cdot\rcf^+,\rho^+)$ is just the extension of
the Lie algebroid $(A^\ast,\lcf\cdot,\cdot\rcf_*,\rho_*)$ by the
vector bundle $M\times\R\to M$ associated with the representation
$\nabla$ and the $2$-cocycle $P_0$. }\end{remark}

Now, one can easily prove the following result.
\begin{proposition}
Let $\pi:A\to M$ be a vector bundle. Consider
$\lcf\cdot,\cdot\rcf^{+}:\Gamma(A_1^*)\times \Gamma(A_1^*)\to
\Gamma(A_1^*)$ (resp. $\rho^{+}:\Gamma(A_1^*)\to C^\infty(M,\R))$
a bracket over $A_1^*$ (resp. a mapping ) given as in (\ref{*}).
Then, $(A_1^*, \lcf\cdot,\cdot\rcf^{+}, \rho^{+})$ is an
almost-special Lie algebroid if and only if the following
sentences are satisfied:
\begin{enumerate}
\item
$(A^*,\lcf\cdot,\cdot\rcf_*,\rho_*)$ is a Lie algebroid,
\item
$X_{0}$ defines a $1$-cocycle of
$(A^*,\lcf\cdot,\cdot\rcf_*,\rho_*),$
\item  $P_0:\Gamma(A^*)\times \Gamma(A^*)\to C^\infty(M,\R)$
is a  skew-symmetric $C^\infty(M,\R)$-bilinear mapping such that
$d_*P_0=-X_{0}\wedge P_0,$ where $d_*$ denotes the differential of
the Lie algebroid $A^*$.
\end{enumerate}
\end{proposition}
Using this result, Corollary \ref{strong-almost} and the local
expression of $\Lambda_A$ and $E_A$, we have that the
strongly-affine Jacobi structure $(\Lambda_A,E_A)$ over $A$
associated with an almost-special Lie algebroid structure over
$A_1^*$ is given by
\begin{equation}\label{LA-EA}
\Lambda_A=\bar\Lambda^*_A -P_0^v + \Delta_A\wedge X_{0}^v,\;\;\;\;
E_A=-X_{0}^v,
\end{equation}
where $\bar\Lambda^*_A$ is the linear Poisson structure over $A$
induced by the Lie algebroid $(A^*,\lcf\cdot,\cdot\rcf_*,\rho_*)$,
$\Delta_A$ is the Liouville vector field of $A$ and $X_{0}^v$
(resp. $P_0^v$ ) is the vertical lift of $X_{0}$ (resp.  $P_0$).

\begin{corollary}\label{1-2}
There is a one-to-one correspondence between strongly-affine
Jacobi structures on a vector bundle $\pi:A\to M$ of rank $n,
n>0,$ and Jacobi algebroid structures
$((\lcf\cdot,\cdot\rcf_*,\rho_*),X_0)$ on $A^*$ with a
skew-symmetric $C^\infty(M,\R)$-bilinear map
$P_0:\Gamma(A^*)\times \Gamma(A^*)\to C^\infty(M,\R)$ such that
$d_*P_0=-X_{0}\wedge P_0,$ where $d_*$ denotes the differential of
$A^*$.
\end{corollary}

\begin{remark}{\rm If $P_0=0$ in Corollary \ref{1-2},
 then we recover the one-to-one
correspondence  between linear-homogeneous Jacobi structures on a
vector bundle $\pi:A\to M$ and Jacobi algebroid structures on
$A^*$ (see Section \ref{Examples}, Example $4$). }
\end{remark}

\section{The characteristic foliation of a strongly-affine Jacobi
structure on a vector space} \setcounter{equation}{0}

Let ${\mathfrak g}$ be a real vector space of finite dimension and
$\bar\Lambda_{\mathfrak g}$ be a linear Poisson structure on
${\mathfrak g}$. Then, ${\mathfrak g}^*$ is a Lie algebra. Denote
by $G^*$ a connected and simply connected Lie group with Lie
algebra ${\mathfrak g}^*$. Then, the leaves of the symplectic
foliation associated with $\bar\Lambda_{\mathfrak g}$ are the
orbits of the coadjoint representation associated with $G^*$. In
this section we will obtain the corresponding result  in the
Jacobi setting.

First  of all, we must replace the terms linear and Poisson by the
terms affine and Jacobi, respectively. So, suppose that we have an
affine Jacobi structure  $(\Lambda_{\mathfrak g},E_{\mathfrak g})$
on ${\mathfrak g}.$ Then, $({\mathfrak g}^+=Aff({\mathfrak
g},\R)={\mathfrak g}^*\times \R,[\cdot,\cdot]^+)$ is a Lie algebra
such that the mapping
\[
-X^{(\Lambda_{\mathfrak    g},E_{\mathfrak    g})}:{\mathfrak
g}^+\to {\mathfrak X}({\mathfrak   g}),\;\;\;\;\;   a    \mapsto
-X^{(\Lambda_{\mathfrak g},E_{\mathfrak
g})}(a)=-X^{(\Lambda_{\mathfrak g},E_{\mathfrak g})}_a
\]
is a Lie algebra antihomomorphism. Denote by $G^+$ a connected and
simply connected Lie group with Lie algebra $({\mathfrak
g}^+,[\cdot,\cdot]^+)$. In general, there does not exist  a global
action of $G^+$ on ${\mathfrak g}$ whose associated infinitesimal
action to be $-X^{(\Lambda_{\mathfrak g},E_{\mathfrak g})}.$ In
fact, in general, if $a\in {\mathfrak g}^+$ then
$-X^{(\Lambda_{\mathfrak g},E_{\mathfrak g})}_a$ is not complete
(see, for instance, (\ref{ejemplo})).

If, "additionally", we suppose that $(\Lambda_{\mathfrak
g},E_{\mathfrak g})$ is strongly affine, then
$-X^{(\Lambda_{\mathfrak g},E_{\mathfrak g})}:{\mathfrak g}^+\to
Aff({\mathfrak g},{\mathfrak g})$ defines an affine representation
of ${\mathfrak g}^+$ on ${\mathfrak g}$ in the sense of \cite{LM}.
Therefore, using a result of Palais \cite{P}, one can prove that
there is an affine representation $\overline{\mbox{Coad}}:
G^+\times {\mathfrak g}\to {\mathfrak g},$  such that the
associated affine representation of ${\mathfrak g}^+$ on
${\mathfrak g}$, $\overline{\mbox{coad}}:{\mathfrak    g}^+\times
{\mathfrak     g}\to {\mathfrak g},$ is $-X^{(\Lambda_{\mathfrak
g},E_{\mathfrak g})}.$ Consequently,

\begin{theorem} Let $(\Lambda_{\mathfrak g},E_{\mathfrak g})$ be a
strongly-affine Jacobi structure over a real vector space
${\mathfrak g}$ of finite dimension.   Then, the leaves of the
characteristic foliation associated with the Jacobi structure
${(\Lambda_{\mathfrak g},E_{\mathfrak g})}$ are just the orbits of
the affine representation $\overline{\mbox{Coad}}:G^+\times
{\mathfrak g}\to {\mathfrak g}.$
\end{theorem}

In the following, we will give an explicit description of the Lie
group $G^+$ and the affine representation
$\overline{\mbox{Coad}}$.

Let $(\Lambda_{\mathfrak g},E_{\mathfrak g})$ be a strongly-affine
Jacobi structure on the real vector space  ${\mathfrak g}$. Then,
there exist a Lie algebra structure $[\cdot,\cdot]_*$ over the
dual vector space ${\mathfrak g}^*$,
 a $1$-cocycle $X_0\in {\mathfrak g}$ of $({\mathfrak g}^*,[\cdot,\cdot]_*)$
and a $2$-section $P_0\in \wedge^2 {\mathfrak g}$ such that
\begin{equation}\label{P0}
d_*P_0=-X_0\wedge P_0,
\end{equation}
$d_*$ being the differential  of $({\mathfrak
g}^*,[\cdot,\cdot]_*)$. Moreover, (see (\ref{LA-EA}))
\[
\Lambda_{\mathfrak g}=\bar{\Lambda}_{\mathfrak g} +
\Delta_{\mathfrak g}\wedge X_0^v-P_0^v,\;\;\; E_{\mathfrak
g}=-X_0^v
\]
where $\bar\Lambda_{\mathfrak g}$ is the Lie-Poisson structure on
${\mathfrak g}$ and $\Delta_{\mathfrak g}$ is the radial vector
field on ${\mathfrak g}$.

 Note that $X_0^v$
(resp.  $P_0^v$) is the constant vector field $C_{P_0}$ (resp. the
constant $2$-vector $C_{P_0}$) over ${\mathfrak g}$ defined by
$X_0\in {\mathfrak g}$ (resp. $P_0\in \wedge^2{\mathfrak g}$).

Let $G^*$ be a connected and simply connected Lie group with Lie
algebra ${\mathfrak g}^*$. Since $d_*X_0=0$, then there is a
unique multiplicative function $\sigma_0:G^*\to \R$ such that
\begin{equation}\label{6.1'}
d\sigma_0(e)=X_0, \end{equation}
 $e$ being the identity element of $G^*$. We recall that $\sigma_0:G^*\to \R$
 is multiplicative if $
\sigma_0(gh)=\sigma_0(g) + \sigma_0(h),$ for $g,h\in G^*.$

On the other hand, $P_0\in \wedge^2{\mathfrak g}$ is a $2$-cocycle
in ${\mathfrak g}^*$ with respect to the representation $R_{X_0}$
of ${\mathfrak g}^*$ on $\R$ defined by
\[
R_{X_0}:{\mathfrak g}^*\times \R \to \R,\;\;\;\; \alpha\in
{\mathfrak g}^*\mapsto
R_{X_0}(\alpha)(\lambda)=\lambda\alpha(X_0), \makebox[1cm]{}
\forall \lambda\in \R. \]

In fact, the cohomology complex associated with this
representation is defined as follows. The space of $k$-cochains
$C^k({\mathfrak g}^*,\R)$ consists of skew-symmetric $k$-linear
mappings $P:{\mathfrak g}^*\times \dots^{(k}\dots\times{\mathfrak
g}^*\to \R$ and the cohomology operator $d_{*X_0}$ is given by
\[
\begin{array}{rcl}
(d_{*X_0} P)(\alpha_1,\dots ,\alpha_{k+1})&=&
\displaystyle\sum_{i=1}^{k+1} (-1)^{i+1}
R_{X_0}(\alpha_i)(P(\alpha_1,\dots ,\widehat{\alpha_i},\dots ,
\alpha_{k+1}))\\&&+\displaystyle\sum_{i<j}(-1)^{i+j}
P([\alpha_i,\alpha_j]_*,\alpha_1,\dots ,\widehat{\alpha_i},\dots,
\widehat{\alpha_j},\dots ,\alpha_{k+1}), \end{array} \] for all
$\alpha_i\in {\mathfrak g}^*.$ Then, $d_{*X_0}P=d_*P + X_0\wedge
P$ and, in particular,  we have that $d_{*X_0}P_0=0$ (see
(\ref{P0})).

Denote by ${Z_{R_{X_0}}^k({\mathfrak g}^*,\R)}$ (resp.
${B_{R_{X_0}}^k({\mathfrak g}^*,\R)}$) the space of $k$-cocycles
(resp. $k$-coboundaries) of the above complex and by
$H_{R_{X_0}}^k({\mathfrak g}^*,\R)$ the corresponding cohomology
group.

On the other hand, $R_{X_0}$ is just  the infinitesimal
representation
 associated with the linear representation  $\phi_{\sigma_0}$ of $G^*$ on
 $\R$ given by
\[
\phi_{\sigma_0}:G^*\times \R\to \R\;\;\;\; (g,t)\mapsto
te^{\sigma_0(g)}. \]
 Now, we consider the cohomology complex over
$G^*$ associated with $\phi_{\sigma_0}$. We recall that the space
of $k$-cochains in this complex is the set $C^k(G^*,\R)$ of
differentiable mappings $\varphi:G^*\times \dots^{(k}\dots \times
G^*\to \R.$ The cohomology operator  is given by
\[
\begin{array}{rcl}
\partial_{\sigma_0}\varphi(g_1,g_2,\dots
,g_{k+1})&=&e^{\sigma_0(g_1)}\varphi(g_2,\dots ,g_{k+1})
+\displaystyle\sum_{j=1}^k(-1)^j\varphi(g_1,\dots
,g_{j-1},g_j.g_{j+1},\\ && g_{j+2},\dots ,g_{k+1}) +
(-1)^{k+1}\varphi(g_1,\dots ,g_k),
\end{array}
\]
for $\varphi\in C^k(G^*,\R)$ and $g_1,\dots ,g_{k+1}\in G^*.$
Denote by $H_{\phi_{\sigma_0}}^k(G^*,\R)=
\displaystyle\frac{Z^k_{\phi_{\sigma_0}}(G^*,\R)}{B^k_{\phi_{\sigma_0}}
(G^*,\R)}$ the cohomology $k$-group of $G$ associated with
$\phi_{\sigma_0}.$

Then, there is an isomorphism between $H^2_{R_{X_0}}({\mathfrak
g}^*,\R)$ and $H^2_{\phi_{\sigma_0}}(G^*,\R)$ (see, for instance,
\cite{AI}). In fact, this isomorphism is induced by the
correspondence
\[
\Phi:Z^2_{\phi_{\sigma_0}}(G^*,\R)\to Z_{R_{X_0}}^2({\mathfrak
g}^*,\R),\;\;
\]
where
\begin{equation}\label{relation}
\Phi(\varphi)(\xi,\eta)=\frac{d}{ds}_{|s=0}\frac{d}{dt}_{|t=0}(\varphi(exp(t\xi),
exp(s\eta))-\varphi(exp(s\eta),exp(t\xi)))
\end{equation}
$exp:{\mathfrak g}^*\to G^*$ being the exponential associated with
$G^*.$

Therefore, for $P_0\in Z_{R_{X_0}}^2({\mathfrak g}^*,\R),$ one may
consider $\varphi_0\in Z_{\phi_{\sigma_0}}^2({\mathfrak g}^*,\R)$
such that $P_0$ and $\varphi_0$ satisfy (\ref{relation}).

Using (\ref{6.1'}) and (\ref{relation}), we prove the following
result.
\begin{theorem}
The Lie group $G^+$ is isomorphic to the product $G^*\times \R$
and the multiplication in $G^+=G^*\times \R$ is given by
\[
\begin{array}{rcl}
(g_1,t_1)(g_2,t_2)&=&(g_1g_2,t_1+e^{\sigma_0(g_1)}t_2-\varphi_0(g_1,g_2)),
\end{array}\] for all $(g_1,t_1),(g_2,t_2)\in G^+=G^*\times \R.$
In particular,
\begin{enumerate}
\item   If   $\varphi_0=0$    (that    is,    the    Jacobi    structure
$(\Lambda_{\mathfrak g},E_{\mathfrak g})$ is linear), $G^+$ is the
semi-direct product $G^*\times_{\phi_{\sigma_0}}\R.$
\item
If $\sigma_0=0$ (that is , the Jacobi structure
$(\Lambda_{\mathfrak g},E_{\mathfrak g})$ is Poisson), $G^+$ is a
central extension of $G^*$.
\end{enumerate}
\end{theorem}

Next, we will describe the affine representation
$\overline{\mbox{Coad}}:G^+\times {\mathfrak g}\to {\mathfrak g}.$

\begin{theorem}
Let ${\mathfrak g}$ be a real vector space of finite dimension and
$(\Lambda_{\mathfrak g},E_{\mathfrak g})$ be a strongly-affine
Jacobi structure over ${\mathfrak g}$. Then, the leaves of the
characteristic foliation associated with $(\Lambda_{\mathfrak
g},E_{\mathfrak g})$ are the orbits of the affine representation
$\overline{\mbox{Coad}}:G^+ \times{\mathfrak g}\to {\mathfrak g}$
given by
\begin{equation}\label{coad}
\overline{\mbox{Coad}}_{(g,t)}(X)=e^{\sigma_0(g)}\mbox{Coad}^{G^*}_gX
+ t X_0-e^{\sigma_0(g)}df_0^{g^{-1}}(e),
\end{equation}
for $(g,t)\in G^+=G^*\times \R$ and $X\in {\mathfrak g}$,  where
$\mbox{Coad}^{G^*}:G^*\times {\mathfrak g}\to {\mathfrak g}$ is
the coadjoint representation associated with $G^*$ and for each
$h\in G^*,$  $f_0^h:G^*\to \R$ is the real function defined by
$f_0^h(g')=\varphi_0(h,g')+\varphi_0(hg',h^{-1}),$ for all $g'\in
G.$
\end{theorem}

{\bf Proof.-} Consider $(\Lambda_{\hat{\mathfrak g}},
E_{\hat{\mathfrak g}})$ the Jacobi structure on $\hat{\mathfrak
g}=({\mathfrak g}^+)^*={\mathfrak g}\times \R$ given by (see
(\ref{r5}))
\[
\Lambda_{\hat{\mathfrak g}}=\bar\Lambda_{\hat{\mathfrak
g}}-\Delta_{\hat{\mathfrak g}}\wedge X_{\zi_{(0,1)}}
^{\bar{\Lambda}_{\hat{\mathfrak g }}}, \;\;\;\; E_{\hat{\mathfrak
g}}=X_{\zi_{(0,1)}}^{\bar{\Lambda}_{\hat{\mathfrak g }}},
\]
where $\bar\Lambda_{\hat{\mathfrak g}}$ is the linear Poisson
structure on $\hat{\mathfrak g}$ associated with the Lie algebra
$({\mathfrak g}^+,[\cdot,\cdot]^+)$ and
$X_{\zi_{(0,1)}}^{\bar\Lambda_{\hat{\mathfrak g}}}$ is the
hamiltonian vector field with respect to
$\bar\Lambda_{\hat{\mathfrak g}}$ associated with the linear
function $\zi_{(0,1)}:\hat{\mathfrak g}\to \R\in {\mathfrak g}^+$
defined by  $\zi_{(0,1)}=(0,1)\in {\mathfrak g}^*\times
\R={\mathfrak g}^+.$
 Denote by $\bar{a}:\hat{\mathfrak g}\to \R$
the linear function associated with  $a=(\alpha,\lambda)\in
{\mathfrak g}^*\times \R={\mathfrak g}^+$. Then, the hamiltonian
vector field with respect to $(\Lambda_{\hat{\mathfrak g}},
E_{\hat{\mathfrak g}})$ associated with $\zi_a$  is
\[
X_{\zi_a}^{(\Lambda_{\hat{\mathfrak g}}, E_{\hat{\mathfrak
g}})}=-a_{\hat{\mathfrak g}}^{\mbox{Coad}^{G^+}} +
X_{\zi_{(0,1)}}^{\bar\Lambda_{\hat{\mathfrak g}}}
(\zi_{a})\Delta_{\hat{\mathfrak g}},
\]
$a_{\hat{\mathfrak g}}^{\mbox{Coad}^{G^+}}$ being the
infinitesimal generator of $a$ associated with the coadjoint
representation $\mbox{Coad}^{G^+}$ of $G^+.$ Consequently, since
${\mathfrak g}=\zi_{(0,1)}^{-1}(1)$, $\Lambda_{\mathfrak
g}=\Lambda_{\hat{\mathfrak g}|{\mathfrak g}}$ and $E_{\mathfrak
g}=E_{\hat{\mathfrak g}|{\mathfrak g}},$ we have that
\begin{equation}\label{restriction}
X_{a}^{(\Lambda_{{\mathfrak g}}, E_{{\mathfrak
g}})}=(X_{\zi_a}^{(\Lambda_{\hat{\mathfrak g}}, E_{\hat{\mathfrak
g}})})_{|{\mathfrak g}}=\big(-a_{\hat{\mathfrak
g}}^{\mbox{Coad}^{G^+}} +
X_{\zi_{(0,1)}}^{\bar\Lambda_{\hat{\mathfrak g}}}
(\zi_a)\Delta_{\hat{\mathfrak g}}\big)_{|\mathfrak g}.
\end{equation}
On the other hand, the coadjoint representation
${\mbox{Coad}}^{G^+}:G^+\times \hat{\mathfrak g}\to \hat{\mathfrak
g}$ is given by
\begin{equation}\label{accionG+}
\mbox{Coad}^{G^+}_{(g,t)}(X,\lambda)=(\mbox{Coad}^{G^*}_g X +
\lambda t e^{-\sigma _0(g)}X_0-\lambda df^{g^{-1}}_0(e), \lambda
e^{-\sigma _0(g)}),
\end{equation}
for all $(g,t)\in G^+=G^*\times \R$ and $(X,\lambda)\in
\hat{\mathfrak g}={\mathfrak g}\times \R.$

From (\ref{restriction}) and (\ref{accionG+}), we obtain that the
action $\widetilde{\mbox{Coad}}:G^+\times \hat{\mathfrak g}\to
\hat{\mathfrak g}$ of $G^+$ on $\hat{\mathfrak g}$ defined by
\begin{equation}\label{accionG+mod}
\widetilde{\mbox{Coad}}_{(g,t)}(X,\lambda)=
e^{\sigma_0(g)}\mbox{Coad}_{(g,t)}^{G^+}(X,\lambda)
\end{equation}
satisfies \[ \widetilde{\mbox{Coad}}_{(g,t)}({\mathfrak
g})={\mathfrak g}\;\;\; \mbox{ and  }\;\;\; (a_{\hat{\mathfrak
g}}^{\widetilde{\mbox{Coad}}})_{|{\mathfrak g}}=
-X_a^{(\Lambda_{\mathfrak g},E_{\mathfrak    g})},\] for all
$a\in{\mathfrak g}^+$ and $(g,t)\in G^+,$
 $a_{\hat{\mathfrak
g}}^{\widetilde{\mbox{Coad}}}$ being the infinitesimal generator
of $a$ associated with the action $\widetilde{\mbox{Coad}}.$
Consequently, the restriction to ${\mathfrak g}$ of
$\widetilde{\mbox{Coad}}_{(g,t)}$ is just
$\overline{\mbox{Coad}}_{(g,t)},$ for all $(g,t)\in G^+.$ Using
this fact, (\ref{accionG+}) and (\ref{accionG+mod}), we obtain
(\ref{coad}).

\hfill$\Box$

Finally, using the above Theorem, we describe the Jacobi structure
on the leaves of the characteristic foliation of a strongly-affine
Jacobi structure over a vector space.

\begin{theorem}
Let ${\mathfrak g}$ be a real vector space of finite dimension and
$(\Lambda_{\mathfrak g},E_{\mathfrak g})$ be a strongly-affine
Jacobi structure over ${\mathfrak g}$. Consider $X\in {\mathfrak
g}$ and $L_{X}$ the leaf of the characteristic foliation over the
point $X$ associated with $(\Lambda_{\mathfrak g},E_{\mathfrak
g})$.
\begin{enumerate}
\item If $E_{\mathfrak g}(X)\notin \#_{\Lambda_{\mathfrak
g}}(T_X^*{\mathfrak g})$ and  $Y\in L_X$ then $$ T_Y
L_{X}=<\{\alpha_Y=\alpha_{\mathfrak g}^{\mbox{Coad}^{G^*}}(Y) +
\alpha(X_0)\Delta_{\mathfrak g}(Y) + (i(\alpha)P_0)^v(Y)\}_
{\alpha\in {\mathfrak g}^*},X_0^v(Y)>$$  and $(\Lambda_{\mathfrak
g},E_{\mathfrak g})$ induces a contact structure $\eta^{L_{X}}$ on
$L_{X}$ given by
\[
\begin{array}{l}
\eta^{L_{X}}(Y)(\alpha_Y)=-\alpha(Y),\;\;\;\;\;
\eta^{L_{X}}(Y)(X_0^v(Y))=-1,
\end{array}
\]
for all $\alpha\in {\mathfrak g}^*.$

\item
If $E_{\mathfrak g}(X)\in \#_{\Lambda_{\mathfrak
g}}(T_X^*{\mathfrak g})$ and $Y\in L_X$ then
$$T_YL_{X}=<\{\alpha_Y=\alpha_{\mathfrak g}^{\mbox{Coad}^{G^*}}(Y)
+ \alpha(X_0)\Delta_{\mathfrak g}(Y) +
(i(\alpha)P_0)^v(Y)\}_{\alpha\in {\mathfrak g}^*}>$$ and
$(\Lambda_{\mathfrak g},E_{\mathfrak g})$ induces a l.c.s.
structure $(\Omega^{L_{X}},\omega^{L_{X}})$ on $L_{X}$ defined by
\[
\begin{array}{l}
\Omega^{L_{X}}(Y)(\alpha_Y, \beta_Y)= [\alpha,\beta]_*(Y)
-P_0(\alpha,\beta),\;\;\;
\omega^{L_{X}}(Y)(\alpha_Y)=-\alpha(X_0),
\end{array}
\] for all $\alpha,\beta\in {\mathfrak g}^*.$
\end{enumerate}
\end{theorem}

{\sc Acknowledgments.} Research partially supported by the Polish
Ministry of Scientific Research and Information Technology under
the grant No. 2 P03A 036 25 and DGICYT grants BFM2000-0808 and
BFM2003-01319. D. Iglesias wishes to thank the Spanish Ministry of
Education and Culture for a FPU grant.

\end{document}